\documentclass[12pt]{article}
\usepackage{latexsym,amsfonts,amsthm,amsmath,amscd,amssymb,color}
\usepackage[dvips]{graphicx}
\usepackage[all]{xy}

\newtheorem{lemma}{Lemma}[section]

\newtheorem{rem}[lemma]{Remark}
\newtheorem{prop}[lemma]{Proposition}

\newtheorem{defn}[lemma]{Definition}

\newcommand{\Ad}{\operatorname{Ad}}

\def\fidi{\hskip5pt \vrule height4pt width4pt depth0pt \par}

\def\sphere{{\mathbb S}}

\def\C{{\mathbb C}}
\def\R{{\mathbb R}}

\def\Z{{\mathbb Z}}
\def\N{{\mathbb N}}

\def\G{{\cal G}}
\def\To{{\mathbb T}}

\def\u{{\mathfrak{u}}}

\def\bbs{\mathbb{S}}

\def\fidi{\hskip5pt \vrule height4pt width4pt depth0pt \par}

\def\gs{\mathfrak{g}^*}

\def\t{\tilde}
\def\t{\mathfrak{t}}
\def\T{{\cal T}}

\newcommand{\BB}{\mathbb}

\newcommand{\FR}{\mathfrak}

\def\R{{\mathbb R}}

\newcommand{\bea}{\begin{eqnarray}}
\newcommand{\eea}{\end{eqnarray}}
\newcommand{\nn}{\nonumber}
\newcommand{\Tr}{\textrm{Tr}}

\newcommand{\sbullet}{\textrm{\tiny{\textbullet}}}

\newcommand{\bra}{\langle}
\newcommand{\ket}{\rangle}

\newcommand{\im}{\textrm{Im}\,}

\def\ga{\alpha}

\def\gc{\gamma}

\def\ep{\epsilon}
\def\gt{\theta}

\def\gs{\sigma}

\def\gl{\lambda}

\def\go{\omega}

\DeclareMathAlphabet{\mathpzc}{OT1}{pzc}{m}{it}

\newtheorem{theorem}{Theorem}[section]
\newtheorem{proposition}[theorem]{Proposition}

\begin{document}

\title{Quantization of Poisson manifolds from the integrability of the modular function}

\author{F. Bonechi\footnote{\small INFN Sezione di Firenze, email: bonechi@fi.infn.it} ,
N. Ciccoli\footnote{\small Dipartimento di Matematica, Univ. di Perugia, email: ciccoli@dipmat.unipg.it} ,
J. Qiu \footnote{\small Math\'ematiques, Universit\'e du Luxembourg, email: jian.qiu@uni.lu},
M. Tarlini\footnote{\small INFN Sezione di Firenze,  email: marco.tarlini@fi.infn.it}
}

\date{June 18, 2013}

\maketitle

\begin{abstract}
We discuss a framework for quantizing a Poisson manifold via the quantization of its symplectic groupoid,
that combines the tools of geometric quantization with the results of Renault's theory of groupoid $C^*$-algebras.
This setting allows very singular polarizations. In particular we consider the case when the modular function is
{\it multiplicatively integrable}, {\it i.e.} when the space of leaves of the polarization inherits a groupoid structure.
If suitable regularity conditions are satisfied, then one can define the quantum algebra as the convolution algebra of the subgroupoid
of leaves satisfying the Bohr-Sommerfeld conditions.

We apply this procedure to the case of a family of Poisson structures on $\C P_n$, seen as Poisson homogeneous spaces of the standard Poisson-Lie
group $SU(n+1)$. We show that a bihamiltoniam system on $\C P_n$ defines a multiplicative integrable model on the symplectic groupoid; we compute
the Bohr-Sommerfeld groupoid and show that it satisfies the needed properties for applying Renault theory. We recover and extend Sheu's description
of quantum homogeneous spaces as groupoid $C^*$-algebras.
\end{abstract}

\thispagestyle{empty}

\eject

\section{Introduction}

The quantization of a Poisson manifold is well understood in terms of the {\it deformation quantization} of the algebra of smooth functions but is much
more elusive if one looks for an {\it operatorial quantization}. In fact Kontsevich formula \cite{Ko}, together with the globalization procedure \cite{CFT},
gives an explicit formula for the deformation quantization of {\it any} smooth Poisson manifold. Nothing similar exists for the operatorial quantization, in
whatever sense of the word.

The concept of {\it symplectic groupoid} was introduced in the mid eighties by Karasev and Weinstein as a tool for quantization of Poisson manifolds.
Indeed, when the Poisson manifold $(M,\pi)$ is {\it integrable} there exists a unique symplectic manifold $\G(M,\pi)$, endowed with a compatible
groupoid structure, that allows to reconstruct the Poisson tensor on its space of units $M$. The compatibility between the groupoid and symplectic structures 
is expressed by requiring that the graph of the multiplication is a lagrangian submanifold of $\overline{\G(M,\pi)}\times\overline{\G(M,\pi)}\times\G(M,\pi)$, 
where $\overline{\G(M,\pi)}$ has the opposite symplectic structure.
If ${\cal H}$ is the space of states obtained by means of some quantization procedure of $\G(M,\pi)$, then, by applying the semiclassical quantization
dictionary, the graph of the multiplication is quantized by a vector in ${\cal H}^*\otimes{\cal H}^*\otimes{\cal H}$. This vector defines an algebra
structure on $\cal H$, that should be considered as the algebra of quantum observables rather than the space of states.

The basic question that one has to address is what kind of quantization we are looking for. It seems a tall order to seek
a map associating an operator to any classical observable, but rather one should look directly for the non commutative algebra of
quantum observables. In fact, from the point of
view of geometric quantization, such a map needs the quantization of any observable, while we know that this is out of reach even in
favourable cases, since the polarization procedure selects a very narrow class of quantizable functions.
Very recently, E. Hawkins in \cite{Hawkins} proposed a general framework based on the geometric
quantization of the symplectic groupoid where the quantization output is a $C^*$-algebra.  Basically, the idea of \cite{Hawkins} is to
look for a {\it multiplicative polarization} $F\subset T_\C\G(M,\pi)$ that allows for the definition of a convolution product on the space of polarized
sections. The best case is that of a real polarization defined by a Lie groupoid fibration $\G(M,\pi)\rightarrow\G_F(M,\pi)$ over some Lie groupoid
$\G_F(M,\pi)$;
the $C^*$-algebra in this case is the convolution algebra of the groupoid $\G_F(M,\pi)$ of lagrangian leaves. The usual problem of geometric quantization is to 
find such polarizations, in particular real ones very often do not exist.

Our starting observation is that thanks to the groupoid structure we can consider much more general polarizations than those needed
by geometric quantization.
In fact, suppose that, instead of a smooth polarization, we just have a topological quotient of groupoids $\G(M,\pi)\rightarrow \G_F(M,\pi)$ with
lagrangian fibres.
Let us call $\G_F^{bs}(M,\pi)$ the subgroupoid of Bohr-Sommerfeld (BS) leaves, {\it i.e.} those leaves with trivial holonomy. In order to define a
convolution algebra it is enough that $\G_F^{bs}(M,\pi)$ admits a Haar system; this can be fulfilled in cases of polarizations that are very singular and hard to use 
in geometric quantization.
In this case, quantization then follows from Renault theory on {\it groupoid $C^*$-algebras} \cite{Renault} and fits in a transparent way with the
Poisson geometry.

Indeed, given a Poisson manifold $(M,\pi)$, there are two cohomology classes in $H_{LP}^\bullet(M,\pi)$, the Lichnerowicz-Poisson cohomology, that are
relevant for quantization. The Poisson tensor itself defines a
class $[\pi]\in H^2_{LP}(M,\pi)$; when a prequantization of the symplectic groupoid is chosen, there is defined a 2-cocycle of the symplectic groupoid with values in
$\sphere^1$  that we call {\it the
prequantization cocycle}. This 2-cocycle can be $1$ if and only if $[\pi]=0$.
Moreover, for every volume form on $M$ it is defined the {\it modular vector field} that measures
the non invariance of the volume form with respect to hamiltonian transformations. It defines a class in $H^1_{LP}(M,\pi)$,
{\it the modular class}, independent on the choice of the volume form: it vanishes only if there exists an invariant volume form. The modular vector
field is integrated to a 1-cocycle of the symplectic groupoid, that we call the {\it modular function}.

It is then natural to require that for our generalized polarization, both the modular function and the prequantization cocycle should descend to
$\G_F(M,\pi)$; the general theory of groupoid $C^*$-algebras tells us how to define the convolution algebra, twisted by the prequantization cocycle; the modular
function is quantized to an automorphism of
the algebra: {\it the modular operator} in noncommutative geometry. Moreover it determines a quasi invariant measure (equivalently the KMS state);
finally we get an Hilbert algebra structure (compatible scalar product and convolution) which corresponds to the algebra on the spaces of
states that we are looking for.

{\it Integrable models} on the symplectic groupoid are a possible source for such singular polarizations and appear very naturally when the Poisson
structure is not unimodular. We say that the {\it modular function is
multiplicatively integrable} if there exists a maximal family of independent hamiltonians in
involution with it such that the space of contour levels inherits the groupoid structure.

A class of examples where such construction is possible is given by compact Poisson-Lie groups and their homogeneous spaces. They are the
semiclassical limit of quantum groups and their quantum homogeneous spaces. A series of papers by
A.J. Sheu \cite{Sheu1,Sh-sphere} relates the quantum spaces to
groupoid $C^*$-algebras. In certain cases (for instance the standard $\C P_{n,q,0}$ and the odd dimensional spheres $\sphere^{2n+1}_q$)
he showed that the $C^*$-algebra is a groupoid $C^*$-algebra; for $SU_q(n)$ or in the non standard structures $\C P_{n,q,t}$, for $t\in(0,1)$, his
characterization of the $C^*$-algebra is just as a subalgebra of a groupoid $C^*$-algebra. His analysis is based on representation
theory of the quantum algebras.

We conjecture that we can interpret the groupoids appearing in Sheu's papers as BS-leaves of some polarization
of the symplectic groupoid integrating the semiclassical structure.  Since they are in general
non unimodular Poisson manifolds, we conjecture that these BS groupoids appear from the multiplicative integrability of the modular function.

In this paper we are going to show that this is actually true for the whole family of Poisson structures on $\C P_n$, thus recovering Sheu's description
for the standard case and extending it in the non standard one. Moreover, as a byproduct of the construction we get the quantization of a class of
Poisson submanifolds; in particular we recover the description of the odd dimensional spheres.

The simplest example is the so called Podle\`s two sphere and has been already studied in \cite{BCST2} (the integrability issues were not explicit there,
for this point of view see the discussion in \cite{BCT}).

In the general case, the relevant integrable model depends on a bihamiltonian structure existing on all compact hermitian spaces. In fact it was
shown in \cite{KRR} that
on compact hermitian symmetric spaces there exists a couple of compatible Poisson structures $(\pi_0,\pi_s)$ defining
a bihamiltonian system. The first one $\pi_0$,  the {\it Bruhat} or {\it standard} Poisson structure, is the quotient of the compact Poisson-Lie group
by a Poisson subgroup and $\pi_s$ is the inverse of the Kirillov symplectic form,
once that we realize the symmetric space
as a coadjoint orbit. By taking a linear combination we get a {\it Poisson pencil}, {\it i.e.} a family of Poisson
homogeneous spaces $\pi_t=\pi_0+t\pi_s$. For $t=(0,1)$ $\pi_t$ is called {\it non standard}. Since one of two Poisson structures is symplectic, we can
define a maximal set of hamiltonians that are in involution with respect to every Poisson
structure of the pencil (for a short account of this see \cite{Foth}).

We analyze here the case of complex projective spaces and we are interested in the case where $t\in[0,1]$, when $\pi_t$ is degenerate.
The integrable model corresponds to the toric structure of $\C P_n$ given
by the action of the Cartan $\To^n\subset SU(n+1)$ and was studied in connection
with the bihamiltonian structure in \cite{Foth}. The action of the Cartan on $\C P_n$ is hamiltonian in the symplectic case
with momentum map $c:\C P_n\rightarrow {\mathfrak t}^*_n$, ${\mathfrak t}_n= {\rm Lie} \To^n$, and simply Poisson with respect to the Poisson
structure $\pi_t$ for $t\in [0,1]$. This means that the action lifts to an hamiltonian action of $\To^n$ on the symplectic groupoid $\G(\C P_n,\pi_t)$
with momentum map $h:\G(\C P_n,\pi_t)\rightarrow{\mathfrak t}^*_n$. The couple $(c,h)$, where $c$ is pulled back to $\G(\C P_n,\pi_t)$ with the source map,
defines an integrable model. We show in Proposition \ref{multiplicative_integrable_model} that it is multiplicative.
The modular vector field with respect to the Fubini-Study volume form is just the fundamental vector field of a particular element of
${\mathfrak t}_n$ so that the modular function is in involution with $(c,h)$.

Our main result is Theorem \ref{main_theorem} where we give an explicit description of the Bohr-Sommerfeld groupoid $\G_F^{bs}(\C P_n,,\pi_t)$ and
show that it has a unique Haar system. Moreover, in Proposition 5.4 we show that the maximal Poisson submanifolds, including odd spheres, are
quantized to subgroupoids of
this BS-groupoid, showing the covariant nature of the correspondence. We show that the BS-groupoid for $t=0$ and the Poisson submanifold
$(\sphere^{2n-1},\pi_t)\subset(\C P_N,\pi_t)$ are the same appearing in Sheu's description
of $\C P_{n,q,0}$ and $\sphere^{2n-1}_q$; the result for $t=(0,1)$ conjecturally improves his
description for $\C P_{n,q,t}$.

This is the plan of the paper. In Section 2 we give a very brief review of basic facts on integration of Poisson homogeneous spaces and groupoid
$C^*$-algebras. In Section 3 we introduce the definition of multiplicative integrability of the modular function and the general quantization scheme.
In Section 4 we describe the geometry of $(\C P_n,\pi_t)$ as Poisson
homogeneous space of $SU(n+1)$ equipped with the quasitriangular Poisson-Lie structure. In Section 5 we interpret this family of Poisson structures as a
Poisson pencil arising by the bihamiltonian system. In
Section 6 we discuss the integrable model on the symplectic groupoid of $(\C P_n,\pi_t)$; we compute the Bohr-Sommerfeld leaves, discuss the Haar system
and the quasi invariant measure associated to the modular function; we discuss the quantization of maximal Poisson submanifolds. Finally in Section 7 we 
relate our results to Sheu's description of quantum homogeneous spaces in terms of groupoid $C^*$-algebras.

\bigskip\bigskip
\noindent{\bf Notations}. We will denote a groupoid with $\G=(\G,\G_0,l_\G,r_\G,m_\G,\iota_\G,\epsilon_\G)$; $\G_0$ is the space of units,
$l_\G, r_\G:{\cal G}\rightarrow \G_0$ are the source and target maps, $\G_2\subset \G\times\G$ denotes those elements
$(\gamma_1,\gamma_2)\in{\cal G}_{2}$ such that
$r_\G(\gamma_1)=l_\G(\gamma_2)$. Moreover, $m_\G:\G_{2}\rightarrow\G$ is
the multiplication, $\iota_\G:\G\rightarrow\G$ is the inversion and $\epsilon_\G:\G_0\rightarrow\G$ is the embedding
of units. We denote with ${}_x\G=l_\G^{-1}(x)$
and $\G_x=r_\G^{-1}(x)$ for $x\in \G_0$. We say that $\G$ is source simply connected (ssc) if ${}_x\G$ is connected and simply connected for any
$x\in \G_0$.

We denote with $\G_k$ the space of strings of $k$-composable elements of $\G$, with the convention that $\G_1=\G$.
The face maps are $d_i:\G_s\rightarrow \G_{s-1}$, $i=0,\ldots s$, defined for $s>1$ as
\begin{equation}\label{face_maps}
d_i(\gamma_1,\ldots\gamma_s) = \left\{ \begin{array}{ll} (\gamma_2,\ldots\gamma_s)& i=0\cr
(\gamma_1,\ldots\gamma_i\gamma_{i+1}\ldots) & 0<i<s\cr
(\gamma_1,\ldots\gamma_{s-1})& i=s\end{array}\right.
\end{equation}
and for $s=1$ as $d_0(\gamma)=l_\G(\gamma)$, $d_1(\gamma)=r_\G(\gamma)$. The simplicial coboundary operator
$\partial^*:\Omega^k(\G_{s})\rightarrow\Omega^k(\G_{s+1})$ is defined as
$$
\partial^*(\omega) = \sum_{i=0}^s (-)^i d_i^*(\omega)  \;,
$$
and ${\partial^*}^2=0$.
The cohomology of this complex for $k=0$ is the real valued groupoid cohomology; $s$-cocycles are denoted as $Z^s(\G,\R)$.

If $\Lambda$ is a line bundle on $\G$, then $\partial^* \Lambda$ denotes the line bundle $d_0^* \Lambda \otimes d_1^* \Lambda^* \otimes
d_2^*\Lambda$ over $\G_2$, where $d_i:\G_2\rightarrow\G$ denote the face maps defined in (\ref{face_maps}).

If $S\subset\G_0$ then we denote with $\G_S=l_\G^{-1}(S)\cap r_\G^{-1}(S)$ the subgroupoid of $\G$ obtained by restriction to $S$.
If a group $G$ acts on $X$, the {\it action groupoid} $X\rtimes G$ is defined on $X \times G$ with structure maps as
$l(x,g)=x$, $r(x,g)=xg$, $m[(x,g)(xg,g')]=(x,gg')$, $\epsilon(x)=(x,e)$, $\iota(x,g)=(xg,g^{-1})$, for $x\in X$ and $g,g'\in G$.

A Poisson structure is denoted $(M,\pi)$, where $\pi$ is the two tensor defining the Poisson bracket on the smooth manifold $M$. We denote with
$\pi$ also the bundle map $\pi: T^*M\rightarrow TM$. The Lichnerowicz-Poisson differential is denoted with $d_{LP}(X)=[\pi,X]$; its cohomology
with $H_{LP}(M,\pi)$.

\bigskip
\bigskip

\section{Poisson manifolds, symplectic groupoids and Haar systems}

We recall in this Section the background material that we will need later. In particular we introduce basic definitions and
properties of Poisson manifolds and symplectic groupoids, of Poisson-Lie groups and their homogeneous spaces. Finally we discuss Haar systems on
topological groupoids.

\subsection{Poisson manifolds and symplectic groupoids}
We say that a topological groupoid $\G$ is a Lie groupoid if $\G$ and $\G_0$ are smooth manifolds, all maps are smooth and $l_\G,r_\G$ are
surjective submersions. A {\it symplectic groupoid} is a Lie groupoid,  equipped with a symplectic form $\Omega_\G$, such
that the graph of the multiplication is a lagrangian submanifold of $\overline{\G} \times \overline{\G} \times \G$, where
$\overline{\G}$ means $\G$ with the opposite symplectic structure. There exists a unique Poisson structure on $\G_0$ such that
$l_\G$ and $r_\G$ are Poisson and anti-Poisson morphisms respectively; the Poisson manifold  $\G_0$, in such case, is said to be {\it integrable}. 
Given a Poisson manifold $(M,\pi)$ is always possible to define a (ssc) topological groupoid $\G(M,\pi)$ as the quotient of the space of cotangent paths with
respect to cotangent homotopies. A {\it cotangent path} is a bundle map $(X,\eta):T[0,1]\rightarrow T^*M$ satisfying
$ d X + \pi(\eta) = 0$; see \cite{CraFer,CatFel} for the definitions of cotangent homotopies and the groupoid structure maps.
The obstruction for the existence of the smooth structure making $\G(M,\pi)$ to be the (ssc) symplectic groupoid integrating $(M,\pi)$ has
been studied in \cite{CrainicFernandes}.

Let $(M, \pi)$ be a Poisson manifold. A vector field $\chi\in {\rm Vect}(M)$ satisfying $d_{LP}(\chi)=0$ is called a {\it Poisson vector field}.
Two distinguished classes are relevant for what follows.
The first one is the class $[\pi]\in H^2_{LP}(M,\pi)$ defined by the Poisson tensor itself. Let us assume that
$M$ is orientable, and let us choose a volume form $V_M$ on $M$. The modular vector field
$\chi_{V_M}={\rm div}_{V_M} \pi$ is Poisson; its cohomology class $\chi_{V_M}\in H^1_{LP}(M,\pi)$
does not depend on the choice of the volume form and is called the {\it modular class}.

Let $(M,\pi)$ be integrable so that $\G(M,\pi)$ is the ssc symplectic groupoid integrating it. Every Poisson vector field $\chi\in {\rm Vect}(M)$
can be integrated to a unique real valued groupoid 1-cocycle $f_{\chi}\in Z^1(\G(M,\pi),\R)$ with the formula
$$
f_\chi[X,\eta] = \int_0^1 \langle \eta(t),\chi(X(t))\rangle ~dt~~~.
$$

We call $f_{\chi_{V_M}}$ the {\it modular function}.

\subsection{Integration of Poisson homogeneous spaces}\label{integration_of_poisson_homogeneous_spaces}

A {\it Poisson-Lie} group is a Lie group $G$ endowed with a Poisson structure $\pi_G$ such that the group multiplication
$(G\times G,\pi_G\oplus \pi_G)\rightarrow (G,\pi_G)$ is a Poisson map. This means that for each $g_1,g_2\in G$ we have that
$$
\pi_G(g_1g_2) = l_{g_1}{}_* \pi_G(g_2) + r_{g_2}{}_* \pi_G(g_1)  ~.
$$
As a consequence, there
exists the simply connected dual Poisson-Lie group $(G^*,\pi_{G^*})$, that acts infinitesimally on $G$, the so called {\it dressing action}.
We assume that this action is integrated to a group action (G is called {\it complete});
as a consequence also $G^*$ is complete. If $g\in G$ and $\gamma\in G^*$ we denote the left and right action of $G^*$ on $G$ as ${}^\gamma g$ and
$g^\gamma$, respectively; analogously the left and right actions of $G$ on $G^*$ are denoted as ${}^g\gamma$ and
$\gamma^g$, respectively.

Any Poisson--Lie group is integrable and its symplectic groupoid has been first described in \cite{LuWe2}. If $G$ is complete then
$\G(G,\pi_G) = G\times G^*$ with $l(g,\gamma)=g$, $r(g,\gamma)=g^\gamma$, $m(g_1,\gamma_1)(g_1^{\gamma_1},\gamma_2)=(g_1,\gamma_1\gamma_2)$,
$\epsilon(g)=(g,e)$, $\iota(g,\gamma)=(g^\gamma,\gamma^{-1})$.
The non degenerate Poisson tensor reads
\bea
\pi_\G = \pi_G + \pi_{G^*} + \sum_a \ell_{T_a}  r_{T^a},\label{Poisson_double}
\eea
where $\ell_{T_a}$ and $r_{T^a}$ are left and right vector fields on $G$ and on $G^*$ with respect to a choice of dual basis of ${\rm Lie G}$ and ${\rm Lie} G^*$.

The (right) action of $(G,\pi_G)$ on $(M,\pi)$ is Poisson if the action seen as a map $(M\times G,\pi\oplus \pi_G)\rightarrow (M,\pi)$ is Poisson.
This means that for each $m\in M$ and $g\in G$ we have that
$$\pi(mg)=g_*\pi(m)+ m_*\pi_G(g)~~~~~.$$
If the action is transitive then $(M,\pi)$ is called a {\it Poisson homogeneous space}. If there exists $m\in M$ where $\pi(m)=0$ then one can show that the
stability subgroup $H_{m}\subset G$ is coisotropic and viceversa if a closed subgroup $H\subset G$ is coisotropic then there exists a unique Poisson
structure $\pi_{H\backslash G}$ on $H\backslash G$ such that the projection map ${\rm pr}_H$ is Poisson. We call these homogeneous spaces {\it embeddable}.
We know that $({\rm Lie H})^\perp\subset({\rm Lie G})^*={\rm Lie} G^*$ is a Lie subalgebra; let us denote with
$H^\perp \subset G^*$ the connected subgroup integrating $({\rm Lie H})^\perp$, which we assume to be closed. Since $H^\perp$ is coisotropic then
there exists a canonical Poisson structure $\pi_{G^*/ H^\perp}$ on $G^*/H^\perp$; let ${\rm pr_{H^\perp}}$ denote the projection map.
Moreover, the dressing action of $H$ preserves $H^\perp$, and viceversa.

Let us recall here how the integration of embeddable homogeneous spaces is done in \cite{BCST}.
The left $G$ action on itself can be lifted to a symplectic action of $\G(G,\pi_G)$; this action admits a Lu momentum map $J:\G(G,\pi_G)\rightarrow G^*$, $J(g,\gamma)={}^g\gamma$.
If we denote $J_H={\rm pr}_{H^\perp}\circ J$, then the symplectic groupoid integrating $(H\backslash G,\pi_{H\backslash G})$ is the symplectic reduction
$H\backslash J_{H}^{-1}({\rm pr}_{H^\perp}(e))$. The symplectic groupoid structure descends from $\G(G,\pi_G)$. We then get the following description
$$\G(H\backslash G,\pi_{H\backslash G})= \{({\rm pr}_H(g),\gamma), g\in G, ~ {}^g\gamma\in H^\perp\}~~~.$$

Let us denote with $G{}_H\!\!\times_H\! H^\perp$ the fibre bundle associated to the homogeneous principal bundle with the dressing action of $H$ on $H^\perp$.
The map $L:\G(H\backslash G, \pi_{H\backslash G})\rightarrow G{}_H\!\!\times_H\! H^\perp$ defined as $L({\rm pr}_H(g),\gamma)=[g,{}^g\gamma]$ is a diffeomorphism
commuting $l$ with the bundle projection; as a consequence $l$ is a fibration. By taking $L\circ\iota$ we prove the same result for the target map $r$.

\subsection{Haar systems}

Let us recall basic facts from \cite{Renault} (see also the very brief description given in Section 4 of \cite{BCST2}).
Let $\T$ be a topological groupoid and let $C_c(\T)$ denote the space of continuous
functions with compact support. A {\it left Haar system} for $\T$ is a family of measures
$\{{}_x\lambda, x\in \T_0\}$ on $\T$ such that
\begin{itemize}
\item[$i$)] the support of ${}_x\lambda$ is ${}_x\T=l^{-1}_\T(x)$;
\item[$ii$)] for any $f\in C_c(\T)$, and $x\in \T_0$, $\lambda(f)(x)=\int_\T f d{}_x\lambda$, defines
$\lambda(f)\in C_c(\T_0)$;
\item[$iii$)] for any $\gamma\in\T$ and $f\in C_c(\T)$, $\int_\T f(\gamma\gamma')d
{}_{r_\T(\gamma)}\lambda(\gamma')=$ $\int_\T f(\gamma')
d{}_{l_\T(\gamma)}\lambda(\gamma')$.
\end{itemize}

The composition of ${}_x\lambda$ with the inverse map will be denoted as $\lambda_x$; this family
defines a {\it right Haar system} $\lambda^{-1}$.
Any measure $\mu$ on the space of units $\T_0$ induces measures $\nu,\nu^{-1}$ on the whole  $\T$
through
$$
\int_\T f d\nu = \int \lambda(f)d\mu \; ,\qquad\int_\T f d\nu^{-1}=\int \lambda^{-1}(f)d\mu\;.
$$
The measure $\mu$ is said to be {\it quasi-invariant} if $\nu$ and $\nu^{-1}$ are equivalent measures;
in this case the Radon-Nikodym derivative $D_\mu=d\nu_\mu/d\nu_\mu^{-1}$ is called the {\it modular function} of $\mu$.
The function $\log D_\mu\in Z^1(\T,\R)$ turns out to be a groupoid 1-cocycle with values in $\R$ and its
cohomology class depends only on the equivalence class of $\mu$.

Let $\zeta\in Z^2(\T,\sphere^1)$ be a continuous 2-cocycle. One can define the convolution
algebra $C_c(\T,\zeta)$ and correspondingly the groupoid $C^*$-algebra $C^*(\T,\zeta)$. When $\zeta=1$ we denote
the convolution and $C^*$-algebras with $C_c(\T)$ and $C^*(\T)$.
Moreover, given a quasi invariant measure $\mu$ one can define the structure of left Hilbert algebra on $L^2(\T,\nu_\mu^{-1})$ such
that the modular operator is exactly given by multiplication with $D_\mu$.

\smallskip
We are interested in understanding when a given topological groupoid has a (possibly unique) Haar system. We will use the following result.
\begin{defn}
A locally compact groupoid is {\it \'etale} if its unit space is open.
\end{defn}

As a consequence, ${}_x\T$, $\T_x$ are discrete for all $x\in\T_0$ and if there exists a left Haar system then it is equivalent to the counting measure.
We will use in particular the following result (Prop. 2.8 in \cite{Renault}):

\begin{prop}
\label{renault_proposition}
A topological groupoid is {\'etale} and admits a left Haar system if and only if $l_\T$ is a local homeomorphism.
\end{prop}

\bigskip
\bigskip

\section{The quantization framework}

Let $(M,\pi)$ be a Poisson manifold, $\G(M,\pi)$ be its symplectic groupoid, that is assumed to be smooth (this was referred in the previous Section
as integrability of $(M,\pi)$; it must not be confused with the notion of integrability of the modular function that we are going
to discuss in this Section). Let us fix a volume form $V_M$ and let
$f_{\chi_{V_M}}\in C^\infty(\G(M,\pi))$ be the modular function.

We are interested in studying cases where the hamiltonian dynamics of the modular function is integrable, {\it i.e.}
$f_{\chi_{V_M}}$ Poisson commutes with a set $F=\{f_i\}_{i=1}^{\dim\G/2}$  of
hamiltonians $f_i\in C^\infty(\G(M,\pi))$ in involution. Such functions should be generically non degenerate in
the usual sense, {\it i.e.} $df_1\wedge\ldots \wedge df_n\not = 0$ on an open and dense subset of $\G(M,\pi)$.

Let us denote with $\G_F(M,\pi)$ the topological space of connected contour levels of $F$. We remark that if we denote
by ${\cal P}_F=\langle X_{f_i}\rangle\subset T\G(M,\pi)$ the (generalized) real polarization defined by the
hamiltonian vector fields $X_{f_i}$ of the $f_i$'s, we have that $\G_F(M,\pi)$ is not in general $\G(M,\pi)/{\cal P}_F$. This is what happens for instance
in \cite{BCST2}. We will not necessarily require that $\G_F(M,\pi)$ be smooth, we will see later what kind of regularity we need.

\begin{defn}
\label{multiplicative_integrability_modular_function}
The modular function $f_{\chi_{V_M}}$ is said to be multiplicatively integrable if there exists on $\G_F(M,\pi)$ a topological groupoid structure such that the
quotient map $\G(M,\pi)\rightarrow\G_F(M,\pi)$ is a groupoid morphism.
\end{defn}

Let us recall the prequantization of the symplectic groupoid: this is the ordinary prequantization of the symplectic structure with additional
requirements that make it compatible with the groupoid structure. The basic results come from \cite{WeXu}, but the presentation is taken
from \cite{Hawkins}. If $\Lambda$ is a line bundle on $\G$ then with $\partial^*\Lambda$ we denote the line bundle
$d_0^*\Lambda\otimes d_1^*\Lambda^*\otimes d_2^*\Lambda$ on $\G_2$.

\begin{defn} \label{dfn:sympl_gpd_preq}
A prequantization of the symplectic groupoid $\G(M,\pi)$ consists of the triple $(\Lambda, \nabla; \zeta_{\Theta_\G})$
where $(\Lambda, \nabla)$ is a prequantization of $\G(M,\pi)$ as a symplectic manifold and $\zeta_{\Theta_\G}$ is
a section of $\partial^* \Lambda^*$ such that:

\begin{itemize}
 \item[$i$)] $\zeta_{\Theta_\G}$ has norm one and is multiplicative, {\it i.e.} it satisfies the following cocycle  condition: for
$(\gamma_1,\gamma_2,\gamma_3)\in \G_3(M)$
 \begin{eqnarray} \label{cocycle_property}
\partial^*\zeta_{\Theta_\G}(\gamma_1,\gamma_2,\gamma_3)&\equiv& \zeta_{\Theta_\G}(\gamma_1, \gamma_2\gamma_3) \zeta_{\Theta_\G}(\gamma_2,\gamma_3)
\zeta_{\Theta_\G}(\gamma_1, \gamma_2)^{-1}
\zeta_{\Theta_\G}(\gamma_1\gamma_2,\gamma_3)^{-1} \cr &=& 1\; .
 \end{eqnarray}
 \item[$ii$)] $\zeta_{\Theta_\G}$ is covariantly constant, {\it i.e.} if $\Theta_\G$ is a (local) connection form for $\nabla$,
then $\zeta_{\Theta_\G}$ (locally) satisfies
 \begin{equation}
\label{parallel_cocycle}
d\zeta_{\Theta_\G} + (\partial^*\Theta_\G) \zeta_{\Theta_\G}=0\;.
\end{equation}
\end{itemize}
\end{defn}

We call $\zeta_{\Theta_\G}$ the {\it prequantization cocycle} defined by $\Theta_\G$. From Theorem 3.2 in \cite{WeXu}
we get that if $\G(M,\pi)$ is prequantizable as a symplectic manifold and $l_\G$-locally trivial, then there exists a unique groupoid prequantization.

The prequantization cocycle can be chosen to be $1$ if $\Omega_\G$ is multiplicatively exact, {\it i.e.} if there exists a primitive $\Theta_\G$ such
that $\partial^*\Theta_\G=0$. Theorem 4.2 in \cite{Crainic} shows that $\Omega_\G$ is multiplicatively exact if and only if the class of the Poisson bivector 
$[\pi]\in H^2_{LP}(M,\pi)$ is zero.

Finally we say that a leaf $\ell\in\G_F(M,\pi)$ satisfies the {\it Bohr-Sommerfeld} conditions if the holonomy of $\nabla$ along $\ell$ is trivial.
We are able to prove that Bohr-Sommerfeld conditions select a subgroupoid of $\G_F(M,\pi)$ only under the following assumption.

\smallskip
\noindent{\bf Assumption 1}: for each couple of composable leaves $(\ell_1,\ell_2)\in \G_F(M,\pi)_2$ the map
$m_\G:\ell_1\times\ell_2\cap \G(M,\pi)_2\rightarrow \ell_1\ell_2$ induces a surjective map in homology.

\smallskip
This assumption is verified by the integrable model we are going to consider in the following Sections.

\smallskip
\begin{lemma}\label{groupoid_bohr_sommerfeld}
If Assumption 1 is satisfied then the set of Bohr-Sommerfeld leaves $\G_F^{bs}(M,\pi)\subset\G_F(M,\pi)$ inherits the structure of a topological groupoid.
\end{lemma}
{\it Proof}. Let $\gamma$ be a loop in the leaf $\ell=\ell_1\ell_2$ and let $\gamma(t)=\gamma_1(t)\gamma_2(t)$, for
$(\gamma_1(t),\gamma_2(t))\in\ell_1\times\ell_2\cap\G(M,\pi)_2$. Then we have that
$$
\int_{\gamma}\Theta_\G = \int_0^1 dt\ \langle \Theta_\G,m_*(\dot\gamma_1\oplus\dot\gamma_2)\rangle = \int_{(\gamma_1,\gamma_2)}\partial^*\Theta_\G +
\int_{\gamma_1} \Theta + \int_{\gamma_2} \Theta~~.
$$
Since $\partial^*\Theta_{\Theta_\G}$ satisfies (\ref{parallel_cocycle}) then we have that if $\ell_1,\ell_2\in \G_F^{bs}(M,\pi)$ then
$\ell_1\ell_2\in\G_F^{bs}(M,\pi)$. \fidi

\smallskip
We refer to $\G^{bs}_F(M,\pi)$ as the {\it groupoid of Bohr-Sommerfeld leaves}.

\medskip
\noindent{\bf Assumption 2}: the groupoid of Bohr-Sommerfeld leaves $\G_F^{bs}(M,\pi)$ admits a Haar system.
\medskip

If Assumption 1 and 2 are satisfied then we get that the groupoid of Bohr-Sommerfeld leaves $\G_F^{bs}(M,\pi)$ admits a Haar system and inherits the
modular cocycle $f_{\chi_{V_M}}$. Moreover, if the prequantization cocycle $\zeta_{\Theta_\G}$ can be chosen to be $F$-invariant it descends to $\G_F^{bs}$ too.
These data give a quantization of the Poisson structure. Indeed, by applying the results of Renault's theory summarized in the previous Section, we
can define the convolution algebra $C_c(\G_F^{bs},\zeta_{\Theta_G})$. Moreover, we can determine the quasi invariant measure
$\mu_{V_M}$ whose modular cocycle is $f_{\chi_{V_M}}$ and then define the left Hilbert algebra on $L^2(\G_F^{bs},\nu_{\mu_{V_M}}^{-1})$.

In the rest of the paper, we will discuss a concrete example where such a quantization exists.

\section{A family of $SU(n+1)$ covariant Poisson structures on $\C P_n$}\label{Poisson_geometry_complex_projective_spaces}

We introduce in this Section the basic example that we are going to quantize. It is a family $\pi_t$, $t\in[0,1]$, of Poisson structures on the complex
projective space $\C P_n$ that are embeddable homogeneous spaces of the standard Poisson structure on $SU(n+1)$.

We will first recall the standard Poisson structure on $SU(n+1)$ (see \cite{LuWe}) and then introduce a family of coisotropic subgroups $U_t(n)\subset SU(n+1)$.

The standard multiplicative Poisson tensor of $SU(n+1)$ can be described intrinsically  integrating  a $\wedge^2\mathfrak{su}$ valued Lie-algebra
2--cocycle; here, in view of later applications, we describe it in terms of the dressing transformation of $SB(n+1,\C)$ on $SU(n+1)$ .
Let us start with the Iwasawa decomposition
\[
SL(n+1,\mathbb{C})=SU(n+1)\times SB(n+1,\mathbb{C})\,, \qquad SB(n+1,\mathbb{C})= A_{n+1} N_{n+1}\, ,
\]
where $A_{n+1}=\{{\rm diag}(v_1,\ldots,v_{n+1}),\ v_i>0,\ v_1v_2\ldots v_{n+1}=1\}$ and $N_{n+1}$ is the group of upper diagonal complex matrices with 1 on the diagonal.
It is easily shown that ${\rm pr}_{A_{n+1}}: SB(n+1,\mathbb{C}) \rightarrow A_{n+1}$ is a group morphism.

The Lie algebras $\mathfrak{su}(n+1)$ and $\mathfrak{sb}(n+1)$ are dual to each other through the pairing $\langle X,\xi\rangle = {\rm Im}\Tr(X\xi)$,
for $X\in\mathfrak{su}(n+1),~ \xi\in\mathfrak{sb}(n+1,\mathbb{C})$.
We will therefore often identify implicitly $\mathfrak{sb}$ with $\mathfrak{su}^*$. We will denote by $p_{1,2}$ the projections from $\mathfrak{sl}$ to $\mathfrak{su}$ and
$\mathfrak{sb}$ respectively.
As said above, any $d\in SL(n+1,\C)$ can be Iwasawa decomposed into $g\gamma$ with $g\in SU(n+1)$ and $\gamma\in SB(n+1,\C)$.  By taking the product $\gamma\cdot g$ and rewriting
it using the Iwasawa decomposition as $\gamma\cdot g={}^{\gamma}g\cdot \gamma^g$, with ${}^{\gamma}g\in SU$ and $\gamma^g\in SB$, one has defined  the \emph{dressing transformation}
\[
SB(n+1,\C)\times SU(n+1))\to SU(n+1)\, ;\quad (\gamma,g)\to{}^{\gamma}g\, .
\]
If $\xi\in\mathfrak{sb}$ then one has $\xi g=gAd_{g^{-1}}\xi$, thus the infinitesimal dressing transformation generated by $\xi$ is
\bea l_g(p_1Ad_{g^{-1}}\xi),\label{vec_dress}\eea
where $g\in SU(n+1)$ and $l_g$ denotes  left multiplication by $g$.
Now regard $\xi$ as in $\mathfrak{su}^*$, one has $Ad_g^*\xi=p_2Ad_{g}\xi$, which can also be seen by using the pairing. The Poisson tensor, when evaluated on two right
invariant 1-forms $(r_{g^{-1}})^*\xi_{1,2}$, is given by the paring of the 1-form $(Ad_g)^*\xi_2$ with the vector field \eqref{vec_dress} generated by $\xi_1$, or in formulae:
\bea \pi\big((r_{g^{-1}})^*\xi_1,(r_{g^{-1}})^*\xi_2\big)= - \bra p_1Ad_{g^{-1}}\xi_1,p_2Ad_{g^{-1}}\xi_2\ket.\label{poisson_K}\eea
From this formula one can show directly that $\pi$ is Poisson and it is furthermore multiplicative $\pi(gh)=(l_g)_*\pi(h)+(r_h)_*\pi(g)$.

Dually, changing the roles of the two groups, one defines the dressing transformation of $SU(n+1)$ on $SB(n+1,\C)$, and can write down a  multiplicative Poisson tensor for the latter. 
In  this way $SB(n+1,\mathbb{C})$ becomes the dual simply connected Poisson-Lie group of $SU(n+1)$.

Let us denote for $t\in[0,1]$
\begin{equation}\label{conjugation_matrix}
\sigma_t = \left( \begin{array}{ccc}\sqrt{1-t} & 0 & \sqrt{t}\cr
                                           0   & {\rm id}_{n-1} & 0\cr
                                            -\sqrt{t}& 0 &\sqrt{1-t} \end{array}\right)\in SU(n+1)~~.
\end{equation}

Using the formula for the Poisson tensor \eqref{poisson_K}, it can be shown that the subgroups
\[
U_t(n)\equiv\sigma_t S(U(1)\times U(n)) \sigma_t^{-1}\subset SU(n+1)
\]
 are coisotropic for all $t\in[0,1]$ (this was first proved in \cite{Sheu2} at the infinitesimal level). Let us denote with $\mathfrak{u}_t(n)={\rm Lie}\ U_t(n)$. 
 Since $U_t(n)$ is coisotropic, $\mathfrak{u}_t(n)^\perp$ is a Lie subalgebra and we will denote with $U_t(n)^\perp\subset SB(n+1,\mathbb{C})$
the subgroup integrating $\mathfrak{u}_t(n)^\perp\subset\mathfrak{sb}(n+1,\BB{C})$. An explicit expression for  $\mathfrak{u}_t(n)^\perp$ can be deduced from \cite{CS}, page 16.
As a consequence of coisotropy there is a uniquely defined Poisson structure $\pi_t$ on $\mathbb{C} P_n=U_t(n)\backslash SU(n+1)$ such that the
quotient map $p_t: SU(n+1)\rightarrow \mathbb{C} P_n$ is Poisson.

\smallskip

It has to be remarked that Poisson structures corresponding to different values of $t$ are not totally unrelated. Let $X_j$, $j=1\ldots n+1$ denote
the homogeneous coordinates of $\mathbb{C} P_n$.

\begin{lemma}\label{tinminust}
The diffeomorphism $\psi: \C P_n\rightarrow \C P_n$ defined as
$$
\psi[X_1,\ldots, X_{n+1}] = [X_{n+1},\ldots, X_1]
$$
sends $\pi_t$ to $-\pi_{1-t}$.
\end{lemma}
{\it Proof}.
Let, in fact, $J=(i\delta_{i,n+1-i})\in SU(n+1)$ be the matrix with all antidiagonal elements equal to $i$. We have clearly that
$\psi$ is just the right multiplication by $J^{-1}$. Then a direct computation shows that
$p_t(JgJ^{-1})=\psi(p_{1-t}(g))$ for all $g\in SU(n+1)$. On the other hand one has that $A\mapsto JAJ^{-1}$ is an anti--Poisson map, 
which therefore induces an anti-Poisson map on complex projective spaces.\fidi

\smallskip
In general $\pi_t$ and $\pi_{t'}$ for $t'\not= 1-t$ are not Poisson or anti Poisson diffeomorphic. This is shown, for instance,
for $\C P_1$ in \cite{BR}.

The family of covariant Poisson structures $\pi_t$ on $\C P_n$ exhibits a sharply different behaviour, depending on whether $t$ is one of the
limiting values $t=0,1$, which will be called the \emph{standard} or \emph{Bruhat--Poisson structure}, or $t\in(0,1)$, that will be
referred to as the \emph{non standard} case. This is best described by analyzing the corresponding symplectic foliation.
It is an interesting fact that such foliation is essentially determined by the images under the projection of the Poisson subgroups of $SU(n+1)$.
For each $k=1,\ldots, n$, let
\[
G_k=S(U(k)\times U(n+1-k))\subset SU(n+1)\, ; \quad P_k(t)= p_t(G_k)\subset \mathbb{C} P_n\, .
\]
Each $G_k$ is a Poisson subgroup so that $P_k(t)$ is a Poisson submanifold of $\C P_n$. As can be easily deduced from \cite{Sto}, Proposition 2.1,
this exhausts the list of maximal (with respect to inclusion) Poisson subgroups of $SU(n+1)$.  We will explicitly describe the
Poisson submanifolds $P_k(t)$.

When $t=0$, {\it i.e.} in the case of the Bruhat-Poisson structure, Poisson submanifolds $P_k(0)\sim\C P_{k-1}$ are contained one inside the other giving rise to the 
following chain of Poisson embeddings:
\begin{equation}\label{poisson_embeddings}
\{\infty\}=\C P_0 \subset \C P_1 \subset \ldots \C P_n
\end{equation}
where $\C P_k$ corresponds to $X_j=0$ for $j>k$ so that $\infty = [1,0,\ldots,0]$.
The symplectic foliation, then, corresponds to the Bruhat decomposition
\begin{equation}\label{bruhat_decomposition}
\C P_n = \bigcup_{i=1}^{n+1} {\cal S}_{i} \;,
\end{equation}
where each ${\cal S}_{i}$ can be described as ${\cal S}_{i}=\{[X_1,\ldots,X_{i},0,\ldots 0],\ X_{i}\not = 0\}\subset \C P_{i-1}$. Thus we have one contractible symplectic 
leaf in each even dimension, which turns out to be symplectomorphic to standard $\C^i$.
The maximal symplectic leaf (the maximal cell) ${\cal S}_{n+1}$ is open and dense in $\C P_n$. On such $2n$-dimensional open cell ${\cal S}_{n+1}$ it is
possible to define the so--called Lu's coordinates $\{y_i\}_{i=1}^n$ as follows (see Appendix \ref{Review_Lu_coordinates} for
more details). Let $z_i = X_i/X_{n+1}$, $1\leq i\leq n$ be the standard affine coordinates on this cell and
 let us define $y_n= z_n$ and for $i<n$
\begin{equation}\label{lu_coordinates}
y_i = \frac{z_i}{\sqrt{1+\sum_{i<j\leq n} |z_j|^2}} ~~~.
\end{equation}

In terms of such coordinates the Bruhat Poisson structure restricted to
${\cal S}_{n+1}$ reads

\begin{equation}
\label{poisson_tensor_lu_coordinates}
\pi_0|_{{\cal S}_{n+1}} = i \sum_{i=1}^{n} (1+|y_i|^2)\partial_{y_i}\wedge\partial_{\bar y_i} \;.
\end{equation}

Let us now move to the non standard case. Most of this description is contained in \cite{CS}, to which we refer for unproven statements.
By direct computation we have that if  $g\in G_k$
then $p_t(g)=[X_1,\ldots , X_{n+1}]=[(1,0,\ldots,0)\sigma_t^{-1}g]\in P_k(t)$, satisfies
\begin{equation}\label{homogeneous_polynomials}
F_{k,t}(X_1,\ldots X_{n+1})=t\sum_{i=1}^k|X_i|^2 - (1-t) \sum_{i=k+1}^{n+1}|X_i|^2 = 0  ~~.
\end{equation}
Since $F_{k+1,t}-F_{k,t}=|X_{k+1}|^2$ we have
\begin{equation}\label{disegual}
F_{n,t}\ge F_{n-1,t}\ge\ldots\ge F_{1,t}\, .
\end{equation} Let $P_{sing}=\bigcup_{k=1}^n P_{k}(t)$,
 we can easily see that $\C P_n\setminus P_{sing}$ is a disconnected union of $(n+1)$ affine spaces.
Each of the $P_k(t)$ is a union of lower--dimensional symplectic leaves; we recall the following facts from \cite{CS}:
\begin{enumerate}
\item  $P_1(t)$ and $P_n(t)$ are Poisson diffeomorphic to the standard Poisson sphere $\bbs^{2n-1}$, by which we mean the Poisson spheres first
described in \cite{VaSo}, having an $\mathbb S^1$ family of symplectic leaves in each even dimension.
\item if $1<k<n$ then $P_k(t)=\bbs^{2k-1}\times \bbs^{2(n-k)+1}/\bbs^1$ where the two odd-dimensional spheres are endowed with the standard odd
Poisson sphere structure and $\bbs^1$ acts diagonally.
\end{enumerate}
It is relevant for what follows to remark here that as Poisson manifolds $P_k(t)$ do not depend on $t$ (up to Poisson diffeomorphism); the
dependence on $t$ is manifested in the embeddings.
What was not detailed in \cite{CS} is the way in which such Poisson submanifolds intersect along lower dimensional symplectic leaves.
Let us remark the following. First, if $k<h$, then $P_k(t)\cap P_h(t)=\bigcap_{i=k}^h P_i(t)$ and this is a simple consequence of inequalities
\eqref{disegual}. Furthermore, from the analysis of the symplectic foliation of odd-dimensional Poisson spheres one can prove that the singular
part of the symplectic foliation of each $P_k(t)$ is just $P_k(t)\cap P_{k-1}(t)\bigcup P_{k}(t)\cap P_{k+1}(t)$; the limiting cases $k=1,n$ of course
behave differently in that only one of the two terms in the union appears. This can be summarized by saying that the foliation resembles more a
stratification where the strata of symplectic leaves of dimension $2m$ are given by the intersections of $n-m$ consecutive Poisson submanifolds $P_k(t)$;
furthermore each $\bigcap_{i=k}^{k+n-m} P_i(t)$ is equal to the intersection of $P_k(t)$ with the linear subspace $X_{k+1}=\ldots=X_{k+n-m-1}=0$. Thus,
for example, the subset of $0$--dimensional symplectic leaves equals
\[
\bigcap_{i=1}^n P_i(t)=\left\{ [X_1,0,\ldots,0,X_{n+1}]\,\big|\, t|X_1|^2-(1-t)|X_{n+1}|^2=0 \right\}\simeq \bbs^1
\]
In each even dimension lower than $2n$ there appears continuous families of symplectic leaves, in sharp contrast with the standard case.

In low dimensions we can say that non standard $\C P_1$ contains a singular $\bbs^1$ of $0$--dimensional leaves of equation $t|X_1|^2-(1-t)|X_2|^2=0$,
separating two $2$--dimensional symplectic cells. As for non standard $\C P_2$, it contains two copies of $\bbs^3$ separating three affine $4$--dimensional
cells: the two embedded Poisson spheres $\bbs^3$ are made of two distinct $\bbs^1$--families of $2$--dimensional symplectic leaves and are intersecting
along a common circle $\bbs^1$ of $0$--dimensional leaves.
The singular part of the symplectic foliation for $\C P_3$ is summarized in the following diagram, which is, of course, also the starting diagram of
leaves of dimension up to $4$ in each higher dimensional $\C P_n$ (and can be easily extended in any dimension):

\vspace{10pt}

\begin{equation}\label{strataPoisson}
\begin{array}{ccccccccc}
\bbs^5&&&&\frac{\bbs^3\times \bbs^3}{\bbs^1}&&&&\bbs^5\\
&\nwarrow&&\nearrow&&\nwarrow&&\nearrow&\\
&&\bbs^3&&&& \bbs^3 &&\\
&&&\nwarrow&&\nearrow&&& \\
&&&&\bbs^1 &&&&
\end{array}
\end{equation}

We summarize the description of the symplectic foliation of $(\C P_n,\pi_t)$ in the following Proposition.

\begin{proposition}
\label{foliation_pi_t}
For $t\in(0,1)$, $(\C P_n,\pi_t)$ decomposes in the following regular Poisson submanifolds
$$
\C P_n = \bigcup_{r,s\geq0,\ r+s=0}^{n} \C P_n^{(r,s)}
$$
where for $r+s=n$, $\C P_n^{(r,s)}$ is an open connected component of $\C P_n\setminus P_{sing}$, and
for $r+s< n$
\begin{equation}\label{cpn_components}
\C P_n^{(r,s)}=\cap_{j=1}^{n-r-s}P_{r+j}(t)\setminus P_{r}(t)\cup P_{n+1-s}(t)~~~~.
\end{equation}
For $r+s=n$, $(\C P_n^{(r,s)},\pi_t)$ is non degenerate. For $r+s<n$, $(\C P_n^{(r,s)},\pi_t)$ is Poisson diffeomorphic to
${\cal S}_{r,s}^{2(r+s)}\times \sphere^1$, where ${\cal S}_{r,s}^{2(r+s)}$ is a
symplectic manifold of dimension $2(r+s)$.
\end{proposition}

We will need the following fact.

\smallskip
\begin{proposition}\label{contractibility_symplectic_leaves}
All symplectic leaves of $(\C P_n,\pi_t)$, for $t=[0,1]$, are contractible.
\end{proposition}
{\it Proof}.
For $t=0,1$ the symplectic foliation corresponds to the Bruhat decomposition so that the result is straightforward. Let us consider $t\in (0,1)$.
Let us discuss first the open symplectic leaves $\C P_n^{(r,n-r)}$ of dimension $2n$. For $r=0,n$
they are defined by $0<F_{1,t}$ and $F_{n,t}<0$ respectively; it is easy to check that they are open disks.
For $0<r<n$, $\C P_n^{(r,n-r)}$ is characterized by the inequality $F_{r,t}<0<F_{r+1,t}$, that implies $F_{r+1,t}-F_{r,t}=|X_{r+1}|^2\not= 0$.
So let us define for $\lambda\in[0,1]$
$$
\epsilon(\lambda,[X_1,\ldots X_{n+1}])=[\lambda X_1,\ldots,\lambda X_{r},X_{r+1},\lambda X_{r+2},\ldots,\lambda X_{n+1}]~~~~.
$$
It is straightforward to verify that $\epsilon$ preserves the inequalities, {\it i.e.}
$\epsilon:[0,1]\times \C P_n^{(r,n-r)}\rightarrow \C P_n^{(r,n-r)}$,
$\epsilon(1,p)=p$ and $\epsilon(0,p)=[0,\ldots,0,1,0,\ldots ,0]$ for all
$p\in \C P_n^{(r,n-r)}$.

Let us consider now the singular part. $\C P_n^{(r,s)}$ is identified by
\begin{equation}\label{rs_stratum}
F_{r,t}< 0=F_{r+1,t}=F_{r+2,t}=\ldots = F_{n-s,t}< F_{n-s+1,t} ~~~.
\end{equation}
By using the definition of the $F$'s it easy to conclude that $[X_1,\ldots,X_{n+1}]\in \C P_n^{(r,s)}$ iff $X_{r+2}= X_{r+3}=\ldots = X_{n-s}=0$,
$X_{r+1},X_{n-s+1}\not = 0$ and
$$
t \sum_{j=1}^{r+1} |X_j|^2 = (1-t) \sum_{j=n-s+1}^{n+1} |X_j|^2  ~~~~.
$$
Note if $r+s=n-1$ no homogeneous coordinate $X$ is zero, but it is easy to adjust the following construction also to this case.
We can then describe $\C P_n^{(r,s)}=\cup_{\mu\in\sphere^1} {\cal S}_{r,s}^{2(r+s)}(\mu)$ where
\begin{eqnarray*}
{\cal S}_{r,s}^{2(r+s)}(\mu) &=& \{(z,\tilde{z})=[z_1,\ldots,z_{r},z_{r+1}=\rho \mu,0,\ldots,0,1,\tilde{z}_{n-s+2},\ldots,\tilde{z}_{n+1}],  \cr
& & ~~~~~\rho>0,~~~ t \sum_{j=1}^{r} |z_j|^2+ t\rho^2= (1-t)(1+\sum_{j=n-s+2}^{n+1}|\tilde{z}_j|^2)\}
\end{eqnarray*}
is a symplectic leaf parametrized by $\mu\in\sphere^1$.
From the properties of the action of the Cartan subgroup, we conclude that
$z_{r+1}\rightarrow z_{r+1}e^{i\phi}$
is a Poisson diffeomorphism from ${\cal S}_{r,s}^{2(r+s)}(\mu)$ to ${\cal S}_{r,s}^{2(r+s)}(\mu e^{i\phi})$. This proves that $\C P_n^{(r,s)}= {\cal S}_{r,s}^{2(r+s)}(1)\times\sphere^1$
as a Poisson manifold.

The following map
$\epsilon: [0,1]\times {\cal S}_{r,s}^{2(r+s)}(\mu)\rightarrow {\cal S}_{r,s}^{2(r+s)}(\mu)$
$$
\epsilon(\lambda,(z,\tilde{z})) =
(\epsilon_0(z,\lambda)z,\lambda\tilde{z}) ~~~~~~
\epsilon_0(z,\lambda)= \sqrt{\lambda^2 + \frac{(1-t)(1-\lambda^2)}{t \sum_{i}^r|z_i|^2+t \rho^2}} $$
contracts ${\cal S}_{r,s}^{2(r+s)}(\mu)$ to the open disk $\{(z_1,\ldots,z_{r}),~ \sum_{i=1}^{r}|z_i|^2 < (1-t)/t\}$. \fidi

\bigskip

Finally let us describe the (ssc) symplectic groupoid integrating $(\C P_n,\pi_t)$. It can be constructed as a symplectic reduction of the Lu-Weinstein symplectic groupoid 
structure on $SL(n+1,\C)$ integrating the Poisson Lie group ($SU(n+1),\pi$), along the lines described in Subsection \ref{integration_of_poisson_homogeneous_spaces}. 
So that we have
\begin{equation}
\label{symplectic_groupoid}
\G(\C P_n,\pi_t) =\{[g\gamma], g\in SU(n+1), \gamma\in SB(n+1,\C), {}^g\gamma\in U_t(n)^\perp\}~~~~,
\end{equation}
where we denoted with $[g\gamma]$ the class of $g\gamma\in SL(n+1,\C)$ under the quotient map $SL(n+1,\C)\rightarrow U_t(n)\backslash SL(n+1,\C)$.
As a smooth manifold, it can be described as the fibre bundle over $\C P_n$ with fibre $U_t(n)^\perp$ associated to the homogeneous principal bundle
with the dressing action of $U_t(n)$ on $U_t(n)^\perp$. We will denote the symplectic form as $\Omega_\G$. Since the $l$ (or $r$) fibers are diffeomorphic
to $U_t(n)^\perp$ that is contractible, the restriction of the symplectic form $\Omega_\G$ to the $l$-fibres is exact; by applying Corollary 5.3 of
\cite{Crainic} we conclude that $\Omega_\G$ is exact.

The homogenous action of $SU(n+1)$ on $\C P_n$ is Poisson; thus it can be lifted to an action on the symplectic groupoid given by right $SU(n+1)$ multiplication.
This action is symplectic and admits Lu momentum map. We are interested in its restriction to the action of the Cartan subgroup $\To^n$, which is hamiltonian in the usual sense.
Under the identification ${\mathfrak t}_n^*={\rm Lie}A_{n+1}$, the momentum map
$h:\G(\C P_n,\pi_t)\rightarrow {\mathfrak t}_n^*$ is defined for each $[g\gamma]\in\G(\C P_n,\pi_t)$ as
\begin{equation}\label{momentum_map_cartan_groupoid}
h[g\gamma] = \log {\rm pr}_{A_{n+1}}(\gamma)   ~.
\end{equation}
For each $ H\in{\mathfrak t}_n$ we denote the component of $h$ along $H$ as $h_H[g\gamma]=\langle H, h[g\gamma]\rangle =
{\rm Im} {\Tr} (H\log{\rm pr}_{A_{n+1}}\gamma)$.
It is easy to check that $h_H$ is the $\R$-valued groupoid cocycle that lifts the Poisson fundamental vector field realizing the ${\mathfrak t}_n$ action
on $\C P_n$.

\section{The bihamiltonian system on $\C P_n$}\label{bihamiltonian_system}

We introduce here the integrable model relevant for the quantization of the covariant Poisson structures introduced in the previous section.
It is, basically, the toric structure on $\C P_n$ defined by the Cartan subgroup action; we feel however necessary to emphasize its
bihamiltonian structure.

Let $\pi_\lambda$ be the inverse of the Fubini-Study symplectic form $\omega_\lambda$ on $\C P_n$ defined in \eqref{Kirillov_Kostant}. We know from
the results of \cite{KRR}, valid on compact Hermitian symmetric spaces, that the Bruhat-Poisson structure $\pi_0$ and
$\pi_\lambda$ are compatible,
{\it i.e.}
$$[\pi_0,\pi_\lambda]=0~,$$
so that  a Poisson pencil $\pi_0+t\pi_\lambda$
$t\in\R$, of $SU(n+1)$-covariant Poisson structures is defined (see also \cite{Foth} for an alternative proof).
We are going to show that the family of Poisson structures $\pi_t$ that we discussed in Section \ref{Poisson_geometry_complex_projective_spaces}
coincide with this Poisson pencil for $t\in[0,1]$.

\smallskip
\begin{lemma}
\label{poisson_pencil}
For each $t\in[0,1]$, $\pi_t=\pi_0 + t \pi_\lambda$.
\end{lemma}

{\it Proof}.
The proof is quite general and depends on the fact that there exists only one Poisson pencil of $SU(n+1)$--covariant Poisson strcutures on  $\C P_n$.
In fact given the two covariant Poisson structures $\pi_t$ and $\pi_0$ then their difference is a $SU(n+1)$--invariant $2$--tensor. In principle there
is no guarantee that such tensor is Poisson; however, as proven in \cite{Sheu2}, up to scalar multiples, the only $SU(n+1)$ invariant
$2$--tensor on $\C P_n$ is the inverse of the Fubini-Study symplectic form.
\fidi

\smallskip

In Appendix \ref{direct_proof} we give a self contained alternative proof of the above Lemma based on a direct computation; this will fix the scalar
multiplying $\pi_\lambda$ appearing in the above proof.

\smallskip

From now on, we will use the notation $\pi_t=\pi_0+t \pi_\lambda$ for $t\in\R$, keeping in mind that $\pi_t$ is obtained as a quotient by
a coisotropic subgroup only for $t\in[0,1]$.

\smallskip

From the general theory of {\it bihamiltonian systems} (see \cite{MM} for an introduction) we get that the recursion operator
$J_\lambda=\pi_0\circ\omega_\lambda$ has vanishing {\it Nijenhuis torsion}, {\it i.e.} for each couple $(v_1,v_2)$ of vector fields
$$
T(J_\lambda)(v_1,v_2)=[J_\lambda v_1,J_\lambda v_2]-J_\lambda ([J_\lambda v_1,v_2) + [v_1,J_\lambda v_2]-J_\lambda [v_1,v_2] )=0  ~.
$$

Let us define $I_k=\frac{1}{k}\Tr J_\lambda^k, k=1\ldots n$. It is shown in Section 4 of \cite{MM} that these global functions satisfy the so
called {\it Lenard recursion relations}
$$
d I_{k+1} = J_\lambda^* d I_k    ~~;
$$
moreover, they are in involution with respect to both Poisson structures $\pi_0$ and $\pi_\lambda$.
Since the modular vector field is the divergence of the Poisson tensor, we can conclude that the
modular vector field with respect to the Fubini-Study volume form $V_{FS}=\frac{1}{n!}\omega_\lambda^{n}$ is independent on $t$.
Moreover, it is a consequence of Theorem 3.5 of \cite{DaFe} that
this modular vector field is the
$\pi_\lambda$ hamiltonian vector field of $I_1$, {\it i.e.}
$$
\chi_{FS} \equiv {\rm div}_{V_{FS}}\pi_t =  \pi_\lambda (d {\rm Tr} I_1)   ~.
$$

We are going to show that the fundamental vector
fields of the Cartan action solve the spectral problem for the recursion operator.

\smallskip
The homogeneous action of the Cartan subgroup $\To^n\subset SU(n+1)$ is hamiltonian with respect to $\pi_\lambda$ with hamiltonian
$c:\C P_n\rightarrow {\mathfrak t}_n^*$. The component of $c$ along $H\in \t_n$ reads
$$
c_H[g] \equiv \langle c[g],H \rangle = {\rm Tr}(g^{-1}\lambda g H)~,
$$
where $\lambda$ is defined in (\ref{highest_weight_cpn}). Let us introduce the following
basis $\langle H_k, k=1,\ldots n\rangle$ for
${\mathfrak t}_n={\rm Lie}\To^n$, where for each $k$ we have
\begin{equation}
\label{cartan_basis}
H_k = \frac{i}{n+1} \left( \begin{array}{cc} (n+1-k)\ {\rm id}_k & 0 \cr
                                  0 & -k \ {\rm id}_{n+1-k} \end{array} \right)  ~.
\end{equation}

In terms of homogeneous and Lu's coordinates we have
\begin{equation}
\label{c_lu_coordinates}
c_k \equiv c_{H_k} + \frac{k}{n+1} = \frac{\sum_{1\leq j \leq k} |X_j|^2}{\sum_j  |X_j|^2}= 1-\Pi_{i=1}^k \frac{1}{1+|y_i|^2} ~,
\end{equation}
respectively. It is clear that
\begin{equation}
\label{standard_simplex}
0 \leq c_1\leq c_2 \ldots \leq c_n\leq 1  ~~~~
\end{equation}
so that the image of the momentum map is the standard simplex $\Delta_n\subset \R^n$.

\smallskip
The action of the Cartan is hamiltonian with respect to $\pi_\lambda$ but only Poisson with respect to $\pi_t$  for all $t\in [0,1]$, {\it i.e.}
the fundamental vector fields of the action define non trivial classes in the Lichnerowicz-Poisson cohomology. It is
possible to define local hamiltonians, as we are going to show in the following Proposition.

\smallskip
\begin{proposition}\label{action_poisson_pencil}
On each open subset of $\C P_n$ where $c_k\not=1-t$, the function
$$b_k= \log |c_k-1+t|$$
is a local hamiltonian for the fundamental vector field of $H_k$, {\it i.e.} $\sigma_{H_k}=\{b_k,-\}_{\pi_t}$. Moreover,
$$
J_\lambda(\sigma_{H_k}) = (c_k-1) \sigma_{H_k}~~~.
$$
\end{proposition}
{\it Proof}. Let $f\in C^\infty(\C P_n)$; then we compute
$$\{b_k,f\}_{\pi_t}= \frac{\partial b_k}{\partial c_k} \{c_k,f\}_{\pi_t}= \frac{\partial b_k}{\partial c_k} \left(\{c_k,f\}_{\pi_0}+t \sigma_{H_k}(f)\right)\ .$$

If we assume also that $X_{n+1}\not= 0$, we can use Lu coordinates $\{y_i\}$, defined in (\ref{lu_coordinates}), where the Poisson
bivector $\pi_0$ takes the form (\ref{poisson_tensor_lu_coordinates}). This assumption does not spoil the global validity of the result and the calculation is
done in this way because the computation in these coordinates is particularly transparent, while a coordinate free computation is lengthy and unilluminating.

Let us first remark that in these coordinates
$$
\sigma_{H_k} = i \sum_{j=1}^k \left( y_j\partial_{y_j}- \overline{y}_j\partial_{\overline{y}_j} \right)\, .
$$
We then compute
$$
\{b_k,f\}_{\pi_t}=\frac{1}{c_k-1 + t} \left(-\frac{1}{\Pi_{j=1}^k(1+|y_j|^2)}+t\right)\sigma_{H_k}(f) = \sigma_{H_k}(f) \ ,
$$
where in the last passage we used (\ref{c_lu_coordinates}).
Finally, if we denote $b_k^0=b_k|_{t=0}=\log |c_k-1|$ then we compute
\begin{eqnarray*}
J_\lambda(\sigma_{H_k}) = \pi_0 dc_k = \frac{\partial c_k}{\partial b_k^0} \pi_0 d b_k^0 =  \frac{\partial c_k}{\partial b_k^0} \sigma_{H_k} =
(c_k-1)\sigma_{H_k} ~,
\end{eqnarray*}
where in the first equality we used the definition of the recursion operator $J_\lambda=\pi_0\circ\omega_\lambda$. \fidi

\smallskip

We then conclude the following result on the modular vector field and class of $(\C P_n,\pi_t)$.

\begin{lemma}\label{computation_modular_vector_field}
The modular vector field of $(\C P_n,\pi_t)$, $t\in \R$, with respect to the Fubini-Study volume form is
$$
\chi_{FS} = \sum_{i=1}^n \sigma_{H_i}~,~
$$
where $\sigma_{H_i}$ is the fundamental vector field of $H_i\in\t_n$. The modular class is non zero for each $t\in[0,1]$.
\end{lemma}

We note that $\sum_{i=1}^n H_i= \frac{1}{2} \sum_{\alpha>0} \alpha$, the half sum of positive roots.

\bigskip
\bigskip

\section{The quantization of $(\C P_n,\pi_t)$}
We introduce here the multiplicative integrable model on the symplectic groupoid $\G(\C P_n,\pi_t)$, where $\pi_t=\pi_0 + t \pi_\lambda$ for
$t\in [0,1]$; we compute its groupoid of Bohr-Sommerfeld leaves and prove that it gives the quantization of $(\C P_n,\pi_t)$.

\subsection{The multiplicative integrable model}
Let us introduce the integrable system. The action of the Cartan $\To^n \subset SU(n+1)$ on $(\C P_n,\pi_\lambda)$ is hamiltonian with momentum map $c$ defined in 
(\ref{c_lu_coordinates}); similarly the right action on $\G(\C P_n,\pi_t)$ is hamiltonian with momentum map $h$ given in (\ref{momentum_map_cartan_groupoid}).
With $\{c_i\}$ and $\{h_i\}$ we denote the components of the momentum maps along the basis $\langle H_i \rangle$
defined in (\ref{cartan_basis}). We remark that since ${\rm pr}_{A_{n+1}}: SB(n+1,\C)=A_{n+1}N_{n+1}\rightarrow A_{n+1}$ is a group homomorphism
the components of $h$ are groupoid cocycles. In particular, as a consequence of Lemma \ref{computation_modular_vector_field} the modular function with
respect to the Fubini-Study volume form is
\begin{equation}
\label{fubini_study_modular_function}
f_{FS} = \sum_{i=1}^n h_i    ~.
\end{equation}

\smallskip
\begin{lemma}
\label{groupoid_integrable_model}
The $2n$ hamiltonians $F=\{l^*(c_i),h_i\}_{i=1}^n$, ,
are in involution.
\end{lemma}
{\it Proof}. Indeed, the $c_i$'s are in involution
with respect to $\pi_t$ for all $t$ due to the bihamiltonian nature of the integrable model and $l^*$ is a Poisson morphism; moreover for each
$[g\gamma]\in \G(\C P_n,\pi_t)$ we compute that
$$ \{l^*(c_i),h_j\}([g\gamma])= \frac{d\ }{dt} l^*(c_i)(g\gamma e^{tH_j})|_{t=0}=\frac{d\ }{dt} c_i(g e^{tH_j})|_{t=0}=0~,$$
where we used the $\To^n$-invariance of $c$ and the fact that $\mathbb{T}^n$ normalizes $SB(n+1)$. \fidi

\smallskip

It is important to compute the following formula.

\smallskip
\begin{lemma}
\label{cartan_cocycle}
On $[g\gamma]\in\G(\C P_n,\pi_t)$, we have that $h_k[g\gamma]=\det_k\gamma$, where $\det_k$ denotes the determinant of the first $k\times k$ minor.
\end{lemma}
{\it Proof}. We have that
$$
h_k[g\gamma]={\rm Im} {\rm Tr}(H_k\log{\rm pr}_{A_{n+1}}\gamma) = \sum_{i=1}^k \log \gamma_{ii} -\frac{k}{n+1}\sum_{i=1}^{n+1} \log \gamma_{ii}
= \log\det\!{}_k\gamma~,
$$
where we used the fact that $\det\gamma=1$. \fidi

\smallskip
Let us first show that the integrable system is {\it multiplicative}. For each $(c,h)\in \R^n\times\R^n$, let
$L_{ch}=\{[g\gamma]\in \G(\C P_n,\pi_t)\ | \ c(g)=c, h(\gamma)=h\}$ denote the contour level set. The proof that $L_{ch}$
is connected is postponed to Proposition \ref{N-orbits} $ii$).

It is easy to check that the map that sends $(c,h)\in \R^n\times \R^n$ to
\begin{equation}
\label{Rn-action}
r(c,h)_i= 1-t + e^{-h_i}(c_i+t -1)   ~,
\end{equation}
defines a $\R^n$-action. Let $\R^n\rtimes \R^n$ denote the action groupoid and let $\R^n\rtimes\R^n|_{\Delta_n}$ denote
the action groupoid restricted to the standard simplex $\Delta_n\subset\R^n$.

\begin{proposition}\label{multiplicative_integrable_model}
The space of contour level sets $\G_F(\C P_n,\pi_t)$ inherits the structure of topological groupoid. For $t\in (0,1)$, it is equivalent to
the following wide subgroupoid of $\R^n\rtimes\R^n|_{\Delta_n}$,
$$\G_F(\C P_n,\pi_t) = \{(c,h)\in\R^n\rtimes\R^n|_{\Delta_n}~|~ c_i=c_{i+1}=1-t \implies~ h_i=h_{i+1}\}~.$$
If $t\in \{0,1\}$ then it is equivalent to the wide subgroupoid
$$\G_F(\C P_n,\pi_t) = \{(c,h)\in\R^n\rtimes\R^n|_{\Delta_n}~|~ c_i=1-t \implies~ h_i=0\}~.$$
As a consequence, the modular function with respect to the Fubini-Study volume form is multiplicatively integrable.
\end{proposition}
{\it Proof}. Let the space of units $(\G_F)_0$ be the contour level sets of $c$. Let $[g\gamma]\in L_{ch}$; we define $l(L_{ch})=L_c$ and
$r(L_{ch})=L_{c(r[g\gamma])}$. We want to show that this map is well defined
and coincides with (\ref{Rn-action}). If $c_i=1-t$, then also $c_i(r([g\gamma]))=1-t$, since fixing $c_i=1-t$ defines the Poisson submanifold $P_i(t)$,
see (\ref{homogeneous_polynomials}) and (\ref{c_lu_coordinates}). If
$c_k(r(g\gamma))\not=1-t$ then the corresponding generator of the Cartan is hamiltonian with respect to $\pi_t$ on every connected component of
$c_k\not = 1-t$ with hamiltonian $b_k$, where $b_k$ is the action variable defined in Proposition \ref{action_poisson_pencil}. As a consequence, when restricted 
to this component $h_k$ has the following \emph{local} exact form in groupoid cohomology,
{\it i.e.} $h_k = l^*(b_k)-r^*(b_k)$, that is, even though $b_k$ is not globally defined, the difference $l^*(b_k)-r^*(b_k)$ is. We then get
$$h_k[g\gamma]=\log \frac{c_k + t -1}{r(c,h)_k+t-1}~~,$$
{\it i.e.} formula (\ref{Rn-action}). Since $h$ is a groupoid cocycle, the map that sends $L_{ch}$ to $(c,h)\in \R^n\rtimes\R^n|_{\Delta_n}$
respects the groupoid multiplication. We have to characterize its image. If $c_k=1-t$ and $[g\gamma]\in L_{ch}$ then we can choose
$g\in S(U(k)\times U(n+1-k))$ so that $\gamma=g^{-1}\lambda g^\gamma$ for some $\lambda \in U_t(n)^\perp$, see (\ref{symplectic_groupoid}). Since $S(U(k)\times U(n+1-k))$ is
a Poisson subgroup, and so it is invariant with respect to the dressing action, we can write in blocks
$$
g =\left( \begin{array}{cc}g_{11}& 0\cr
                       0& g_{22}
          \end{array}\right) ,~~~~
\lambda = \left(\begin{array}{cc}\lambda_{11}&\lambda_{12}\cr
                                  0 & \lambda_{22}\end{array}\right) , ~~~~
g^\gamma =\left( \begin{array}{cc}g_{11}^\gamma& 0\cr
                       0& g_{22}^\gamma
          \end{array}\right)~~~~.
$$
We can then compute
$$
\exp h_k[g\gamma] = {\det}_k(\gamma) = \det(g_{11}^{-1})\det(\lambda_{11}) \det(g_{11}^\gamma) = \det({\lambda_{11}})={\det}_k(\lambda)~,
$$
since the $\det$ is invariant under the dressing action in $U(k)$. Let us recall that $\u_t(n)^\perp=\Ad^*_{\sigma_t}\u_0(n)^\perp$ (see (\ref{conjugation_matrix})),
so that we can write $\lambda=\exp\Ad^*_{\sigma_t}\xi$ for $\xi\in\u_0(n)^\perp$, {\it i.e.} $\xi_{ij}\not= 0$ only if $i=1, j>1$.
A direct computation shows that
$$
{\rm pr}_{A_{n+1}}(\Ad^*_{\sigma_t}\xi) = \sqrt{t(1-t)} \xi_{1,n+1} {\rm diag}(1,0,\ldots,-1)~,~~
$$
so that we see from [\ref{momentum_map_cartan_groupoid}] that $h_k[g\gamma]=\sqrt{t(1-t)}\xi_{1,n+1}$ and the result easily follows. \fidi

\smallskip
\begin{rem}{\rm
The orbits of $\G_F(\C P_n,\pi_t)$ correspond to the decomposition of $\C P_n$ in regular Poisson strata given in  (\ref{bruhat_decomposition}) for
$t=0,1$ and in Proposition \ref{foliation_pi_t} for $t\in(0,1)$. Moreover, the restrictions on $h$ reflect the fact that when $t\in(0,1)$
on each non maximal stratum the Cartan vector fields which are transversal to the symplectic leaves collapse to a single $\sphere^1$ generator, while for $t=0,1$
they are all zero.}
\end{rem}

\medskip

\subsection{The Bohr-Sommerfeld groupoid $\G_F^{bs}(\C P_n,\pi_t)$}
In order to compute the Bohr-Sommerfeld leaves we need to show that $L_{ch}$ is connected and describe $H_1(L_{ch},\Z)$. It is obvious that
for each $c\in\Delta_n$ $L_c$ is a $\To^n$-orbit.

It is easy to produce 1-cycles in $L_{ch}$. The hamiltonian flux of $h_j, j=1\ldots n$ and  of $l^*(b_j)$ for $c_j\not= 1-t$ define a set of 1-cycles in $L_{ch}$.
Indeed, fix $[g\gamma]\in L_{ch}$;
$\alpha_j(t)= [g\gamma e^{tH_j}]$, for $j=1,\ldots n$, is the flux of $h_j$. Let us compute the flux of $l^*(b_j)$, for $c_j\not=1-t$, and show
that it is periodic.

\begin{lemma}\label{beta_cicles}
The hamiltonian flux of $l^*(b_j)$ passing through $[g\gamma]\in L_{ch}$ is given by
$$
\beta_j(t) = [g e^{B_{H_j}[g]t} \gamma] ~,
$$
where
$$B_{H_j}[g]=H_j+\sum_a \ell_{T_a}(b_j)([g]) T^a\in  {\mathfrak {sl}}(n+1,\C)~~~,$$
for some basis $\{T_a\}$ of ${\mathfrak {su}}(n+1)$ and the dual basis $\{T^a\}$ of $\FR{sb}(n+1)$. Moreover, $\beta_j(t)=\beta_j(t+2\pi)$ for each $t$.
\end{lemma}
{\it Proof}. Let $\beta_j(t)=[g_j(t)\gamma_j(t)]$ be the hamiltonian flux of $l^*(b_j)$. Since $l(\beta_j(t))=[g e^{tH_j}]$ then $[g_j(t)\gamma_j(t)] = [ge^{t H_j}\gamma_j(t)]$.
By using the definition of the Poisson tensor (\ref{Poisson_double}), especially the third term thereof that entangles $G$ and $G^*$, we get that
\bea
\dot\gamma_j(t) = (\sum_a \ell_{T_a}(b_j)[ge^{tH_j}] T^a )\gamma_j(t) = e^{-tH_j}(\sum_a\ell_{T_a}(b_j)[g] T^a) e^{t H_j}\gamma_j(t) ~,\label{flow_equation}
\eea
note that in the second equality we used the equivariance of the momentum map, which can be seen explicitly from
\bea \ell_{T_a}(c_j)([g])={\rm Tr}[g^{-1}\lambda g [T_a,H_j]],\label{momen_equiv}\eea
where $\lambda$ is the weight of the fundamental representation of $SU(n+1)$.

The solution to the equation (\ref{flow_equation}) is
$$
\gamma_j(t) = e^{-tH_j} e^{tB_{H_j}[g]}\gamma~.
$$
We just have to show that $\beta_j$ is periodic. It is a direct computation to obtain a block matrix form for $B_{H_j}[g]$: choose the generators
of the Lie algebra of $N_{n+1}\subset SB(n+1,{\mathbb C})$ as $(E^{ij})_{kl}=\delta_{ik}\delta_{jl}$ and $i E^{ij}$ with $i<j$, and the
dual generators $i (E^{ij}+E^{ji})$ and $(E^{ij}-E^{ji})$. For $k,l<j$ or $k,l>j$ the latter commute with $H_j$, and hence from (\ref{momen_equiv})
the corresponding blocks of the matrix $\ell_{T_a}(c_j)([g])$ vanishes, the blocks that remain are
\begin{equation}
\label{b-flux}
B_{H_j} = \frac{i}{n+1} \left( \begin{array}{cc} (n+1-k)\ {\rm id}_k & A\cr
                                  0 & -k \ {\rm id}_{n+1-k} \end{array} \right)  ~.
\end{equation}
It can be diagonalized to $H_j$ with the similarity transformation by
$$\left(\begin{array}{cc} 1& A/(n+1)\cr 0&1 \end{array}\right)~.$$
The periodicity of $\beta_j$ then follows from the observation that $\exp 2\pi H_j\in\Z_n\subset U_t(n)$. \fidi

\smallskip
We are going to show that $\{\alpha_j,\beta_k\}$ generate $H_1(L_{ch},\Z)$. We will use arguments based on the use of the homotopy long exact sequence
associated to a fibration. We need to introduce some definitions and preliminary results.

Let $x\sim_t y$ denote that $x,y\in \C P_n$ live in the same symplectic leaf of $\pi_t$.
Let $\tilde{L}_{ch}=\{(x,y)\in L_c\times L_{r(c,h)}\ | \ x\sim_t y\}\subset L_c\times L_{r(c,h)}$; we have clearly that
\begin{equation}\label{fibration}
(l\times r): L_{ch} \rightarrow \tilde{L}_{ch}~.
\end{equation}
\begin{lemma}\label{Tc}
Let $x\in L_c$ and let $\To_x=\{s\in \To^n\ | xs \sim_t x\}$,
Then we see that
\begin{itemize}
 \item [$i$)] For each $x,y\in L_c$, $\To_x=\To_y\equiv \To_c$.
 \item [$ii$)] For each $x \in L_c$, $T_x (x\To_c) = \langle \sigma_{H_i}(x)\ , c_i\not = 1-t \rangle$.
 \item [$iii$)] The hamiltonian flux of $l^*(b_j)$, for $c_j\not= 1-t$, described in Lemma \ref{beta_cicles} defines an action of $\To_c$ on $L_{ch}$ that commutes
                with the $\To^n$-action.
 \item[$iv$)] $\To_c\times \To^n$ acts transitively on
$\tilde{L}_{ch}$ as $(x,y)(s,k)=(xsk,yk)$, where $(x,y)\in \tilde{L}_{ch},(s,k)\in \To_c\times\To^n$;
the map  (\ref{fibration}) intertwines the actions.
\end{itemize}
\end{lemma}
{\it Proof}. $i$) Since the action of the Cartan is by Poisson diffeomorphisms, its tangent action commutes with the sharp map, so that if $x\sim_t y$ then
$x\tau\sim_t y\tau$ for any $\tau\in \To^n$. Let $x\in L_c$ and $s\in\To_x$; let $y=x \tau$, for some $\tau\in \To^n$. We have that $x\sim_t xs$ implies
$y\sim_t ys$ so that $s\in \To_y$.

$ii$) Let $x\in L_c$ and let ${\cal S}_x$ be the symplectic leaf through $x$. If $c_i\not=1-t$ then $\sigma_{H_i}$ is locally hamiltonian and clearly preserves the leaf.
Let $c_i=1-t$ and let $\sigma_{H_i}(x)\not=0$
so that $d\phi_i(x)$ is defined and non zero. Clearly $\langle \sigma_{H_i}(x),d\phi_i(x)\rangle=1$. We can choose an open $U_x$ containing $x$ such that the points 
$y$ such that $d\phi_i(y)$ is defined and $c_i(y)\not =1-t$ are dense. For such $y$, $\pi_t(d\phi_i(y))=\partial/\partial_{b_i}=(c_i-1+t)\partial/\partial_{c_i}$
that extends to 0 in $x$. So we have that $d\phi_i(x)\in N^*_x{\cal S}_x$ and $\sigma_{H_i}(x)\not\in T_x{\cal S}_x$.

$iii$) This is just a rephrasing of Lemma \ref{beta_cicles}.

$iv$) Let us first show transitivity. Let $(x,y),(x_0,y_0)\in \tilde{L}_{ch}$; then we know that
$(x,y)=(x_0k_1,y_0 k_2)=(x_0k_1k_2^{-1},y_0)(1,k_2)$ for $k_1,k_2\in \To^n$, so that $x_0 k_1 k_2^{-1}\sim y_0\sim x_0$ so that $k_1k_2^{-1}\in \To_c$
and $(x,y)=(x_0,y_0)(k_1k_2^{-1}, k_2)$. It is a trivial fact that $l\times r$ intertwines the $\To^n$-action; the same property with respect to $\To_c$
follows from the properties of the hamiltonian vector fields of $l^*(b_j)$.
\fidi

\smallskip
For each $[g]\in \C P_n$ let us denote $N_{n+1}[g]=\{n\in N_{n+1}\ | \ {}^gn\in U_t(n)^\perp, [g^n]=[g] \}= g^{-1}U_t(n)^\perp U_t(n) g\cap N_{n+1}$; it is
clearly a subgroup of $N_{n+1}$. Next we shall show that this subgroup is of exponential type (and hence contractible) and that two such subgroups are
isomorphic if the base points live in the same $\To^n$-orbit.

\smallskip
\begin{lemma}
\begin{itemize}
\item[$i$)] $N_{n+1}[g]\subset N_{n+1}$ is connected and then of exponential type.
\item[$ii$)] For $s\in\To^n$, the map $n\rightarrow {}^{s^{-1}}n$ is a group isomorphism from $N_{n+1}[gs]\rightarrow N_{n+1}[g]$.
\end{itemize}
\end{lemma}
{\it Proof}. $i$) We can easily see that $N_{n+1}[g]= {}_{[g]}\G(\C P_n,\pi_t)_{[g]}\cap N_{n+1}$. From Lemma 1 in \cite{CrainicFernandes}, we know that
$$\pi_0({}_{[g]}\G(\C P_n,\pi_t)_{[g]}) = \pi_1({\cal S}_{[g]})  ~,$$
where ${\cal S}_{[g]}$ is the symplectic leaf through $[g]$.
Since in our case the symplectic leaves are all contractible (see Lemma \ref{contractibility_symplectic_leaves}), we conclude that $N_{n+1}[g]$ is the
intersection of two connected subgroups of $SB(n+1,\C)$; since
connected subgroups of an exponential group are of exponential type and are closed under the intersection (see \cite{Sa} pg.9) we get the result.
Point $ii$) is obtained by direct computation. \fidi

\smallskip
Let $N_c$ denote the group bundle over $L_c$ with fibre $N_{n+1}[g]$ over $[g]\in L_c$ and let $p:N_c\rightarrow L_c$ denote the projection;
it acts on $L_{ch}$ with anchor $l:L_{ch}\rightarrow L_c$ and action $a: N_c\ {}_p\!\times_l L_{ch}\rightarrow L_{ch}$ given by
\begin{equation}\label{group_bundle_action}
a(n,[g\gamma]) = [gn\gamma]    ~~~~.
\end{equation}
This action preserves the fibres of (\ref{fibration}).

\smallskip
\begin{proposition}
\label{N-orbits}
\begin{itemize}
 \item[$i$)] The map (\ref{fibration}) defines a fibre bundle and the groupoid action of $N_c$ on $L_{ch}\rightarrow \tilde{L}_{ch}$
\begin{equation*}
\xymatrix{ N_c  \ar[d]_{p}
& L_{ch} \ar[dl]^{l} \ar[d]^{l\times r} \\ L_c & \tilde{L}_{ch} }
\end{equation*}
defines a principal groupoid bundle.
 \item[$ii$)] $L_{ch}$ is connected and $H_1(L_{ch},\Z)$ is generated by $\langle \alpha_j,\beta_k\rangle$.
\end{itemize}
\end{proposition}
{\it Proof.} $i$) Let us show first that we have a groupoid principal bundle according to definition in section 5.7 of \cite{MoMr}.
The map $l\times r: L_{ch}\rightarrow\tilde{L}_{ch}$ is a surjective submersion as a consequence of Lemma \ref{Tc} (iv). Next we have to show that the map
$$
(a,{\rm pr}_2): N_c\ {}_{L_c}\!\times_{L_c} L_{ch}\rightarrow L_{ch}\ {}_{\tilde{L}_{ch}}\!\!\times_{\tilde{L}_{ch}} L_{ch}~,~~~
(n,[g\gamma])\rightarrow ([gn\gamma],[g\gamma]) $$
is a diffeomorphism.
Let $([\tilde{g}\tilde{\gamma}],[g\gamma])\in L_{ch}\ {}_{\tilde{L}_{ch}}\!\!\times_{\tilde{L}_{ch}} L_{ch}$, {\it i.e.} $(l\times r)[g\gamma]=(l \times r)[\tilde{g}\tilde{\gamma}]$. 
Then we see that $[g]=[\tilde{g}]$ and $g^{\tilde\gamma}=h g^\gamma$ for $h\in U_t(n)$. Define $n=\tilde{\gamma}\gamma^{-1}\in N_{n+1}$. Since we also have that
${}^g\gamma,{}^g\tilde{\gamma}\in U_t(n)^\perp$ and recalling that the dressing action of $U_t(n)$ preserves $U_t(n)^\perp$ (and viceversa), one shows that
\begin{eqnarray*}
{}^gn&=& {}^g\tilde{\gamma}\  {}^{g^{\tilde{\gamma}}}\gamma^{-1}= {}^g\tilde{\gamma}\  {}^h({}^{g^{\gamma}}\gamma^{-1})=
{}^g\tilde{\gamma}\  {}^h({}^{g} \gamma)^{-1} \in U_t(n)^\perp ~,\cr
g^n&=& g^{\tilde{\gamma}\gamma^{-1}} = (hg^\gamma)^{\gamma^{-1}}= h^{({}^g\gamma)^{-1}}g =h'g~,
\end{eqnarray*}
for $h,h'\in U_t(n)$. We conclude that $n\in N_{n+1}[g]$ and the map $([\tilde{g}\tilde{\gamma}][g\gamma])\rightarrow (n,[g\gamma])$ defines an inverse of
$(a,{\rm pr}_2)$. The map $\delta_{N_c}:L_{ch}\ {}_{\tilde{L}_{ch}}\!\!\times_{\tilde{L}_{ch}} L_{ch}\rightarrow N_c$,
$([\tilde{g}\tilde{\gamma}][g\gamma])\rightarrow \tilde{\gamma}\gamma^{-1}$, is called
the {\it division map}.

Due to Lemma \ref{Tc} $iv$), $\tilde{L}_{ch}$ is homeomorphic to $\To_c\times\To^n/T$ for $T$ being the stability subgroup. Let us fix
$[g_0\gamma_0]\in L_{ch}$ and let $(x_0,y_0)=(l\times r)[g_0\gamma_0]$ and $V_{0}$ an open containing $(x_0,y_0)$ and admitting a section
$(s,k): V_{0}\rightarrow \To_c\times\To^n$ such that $(x,y)=(x_0,y_0)(s(x,y),k(x,y))$ for each $(x,y)\in V_0$.
Let $[g\gamma]\in L_{ch}|_{V_0}$ and $(x,y)=(l\times r)[g\gamma]$; then $[g\gamma](s(x,y),k(x,y))^{-1}\in (l\times r)^{-1}(x_0,y_0)$ and
$[g\gamma]\rightarrow \delta_{N_c}([g_0,\gamma_0],[g\gamma](s(x,y),k(x,y))^{-1})\in N_{[g_0]}$ gives the local trivialization.

$ii$) We apply the long exact homotopy sequence to the fibre bundle (\ref{fibration}). Since the fibre is contractible we get that
$\pi_k(L_{ch})=\pi_k(\tilde{L}_{ch})$, in particular $L_{ch}$ is connected. The same homotopy sequence applied to
$T\rightarrow\To_c\times\To^n\rightarrow \tilde{L}_{ch}$, since $T$ is connected, shows that
the map $\pi_1(\To_c\times\To^n)\rightarrow \pi_1(\tilde{L}_{ch})=\pi_1(L_{ch})$, sending the generators of $\To_c\times\To^n$ to $\{\alpha_j,\beta_k\}$,
is surjective. The statement on homology easily follows. \fidi

\smallskip
We are now ready for the computation of the groupoid of Bohr-Sommerfeld leaves. From the discussion at the end of Section 3, we know
that the symplectic form $\Omega_\G$ is exact. In the following Lemma we show that there exists a primitive $\Theta_\G$, {\it i.e.} a one form $\Theta_\G$ satisfying
$d\Theta_\G=\Omega_\G$, adapted to the momentum map $h$ of the action of the Cartan $\To^n$. Remark that the fact that $h$ is
an $\R^n$-valued cocycle plays a role.

\smallskip
\begin{lemma}
\label{adapted_potential}
There exists a primitive $\Theta_\G$ of the symplectic form $\Omega_\G$ such that for each
$H\in{\mathfrak t}_n$
$$
\iota_{\sigma_H}\Theta_\G = - h_H ~~~,
$$
where $\sigma_H$ is the fundamental vector field corresponding to $H$.
\end{lemma}
{\it Proof.} Since $\To^n$ is compact, by averaging, we can always find a $\To^n$-invariant primitive $\Theta_\G$,  {\it i.e.} $L_{\sigma_H}\Theta_\G=0$.
By using the definition of the momentum map we easily see that
$$
\iota_{\sigma_H}\Theta_\G+h_H= a_H \in \R~~.
$$
Since the submanifolds of units $\G_0=\C P_n$ is lagrangian then $\Theta_\G|_{\G_0}$ is closed and since $H^1(\C P_n)=0$ then
$\Theta_\G|_{\G_0}$ is exact. The $\To^n$-action on $\G(\C P_n,\pi_t)$ preserves $\G_0$; moreover for any orbit
$\gamma_H\subset \C P_n$ of $\sigma_H$ we have that $h_H|_{\gamma}=0$ because $h_H$, being an $\R$-valued cocycle, is zero on the units. We then conclude that
$$
\int_{\gamma_H} \Theta_\G = a_H = 0   ~~~.
$$
\fidi

\smallskip
From now on, $\Theta_\G$ will denote a fixed primitive satisfying the condition of Lemma \ref{adapted_potential}.

\smallskip
\begin{theorem}
\label{main_theorem}
For $\hbar\in\R$, the groupoid of Bohr-Sommerfeld leaves is the following subgroupoid of $\G_F(\C P_n,\pi_t)$ described in Proposition
\ref{multiplicative_integrable_model}
$$
\G_F^{bs}(\C P_n,\pi_t) =\{(c,h)\in\G_F(\C P_n,\pi_t) \ |\ \log|c_k-1+t|\in \hbar\Z,\ h_k\in \hbar\Z\}    ~~~~~~.
$$

It is an \`etale groupoid and admits a unique left Haar system for all $t\in[0,1]$. The generators of the Cartan subgroup are quantized to the groupoid
cocycles $h_i(c,h)=h_i$; in particular the modular cocycle reads
$$
f_{FS}(c,h)=\sum_{i=1}^n h_i~~\in \hbar \Z.
$$
\end{theorem}

{\it Proof}. Let $L_{ch}\in\G_F(\C P_n,\pi_t)$ and let $c_k\not =1-t$.  Let $\Theta_\G$ be a primitive of $\Omega_\G$ as in Lemma \ref{adapted_potential};
we have to pair it with the generators $\langle\alpha_j,\beta_k\rangle$ of $H_1(L_{ch},\Z)$. We then compute the Bohr-Sommerfeld
conditions as
$$
\int_{\beta_k} \Theta_\G = \int_{\tilde{\beta}_k} \Omega_\G|_{\G_{[g]}}=\int_{\tilde{\beta}_k}l^*\omega_{{\cal S}_{[g]}}=
\int_{l(\tilde{\beta}_k)} \omega_{{\cal S}_{[g]}} = b_k 2\pi= 2\pi\log|c_k-1+t| \in 2\pi\hbar\Z~,
$$
where $\tilde{\beta}_k:{\mathbb D}\rightarrow \G(\C P_n,\pi_t)_{[g]}$ extends $\beta_k$ to the disk and in the second equality we use the fact
that $\Omega|_{\G_{[g]}}=l^*\omega_{{\cal S}_{[g]}}$. Since
$\iota_{\sigma_H} \Theta_\G = - h_H$, for each $H\in {\mathfrak t}_n$, we easily compute
$$
\int_{\alpha_j} \Theta_\G=\int_0^{2\pi}\iota_{\sigma_{H_j}}\Theta_\G\ dt=-2\pi h_j \in 2\pi \hbar\Z~~.
$$
The fact that the set of Bohr-Sommerfeld leaves $\G_F^{bs}(\C P_n,\pi_t)$ is a subgroupoid is a straightforward computation.

Let $(c,h)\in \G_F^{bs}(\C P_n,\pi_t)$ and let $I_{\epsilon}(c,h)=f_{FS}^{-1}(I_\epsilon(f_{FS}(c,h)))$ the open neighborhood of
$f_{FS}(c,h)$ of radius $\epsilon>0$. It is clear that for $\epsilon < 1$,
$I_{\epsilon}(c,h)=\{(c',h)\in\G_F^{bs}(\C P_n,\pi_t)\}$ so that $l$ is a local homeomorphism. By applying Proposition
\ref{renault_proposition} we conclude that $\G_F^{bs}(\C P_n,\pi_t)$ is \'etale and admits a unique left Haar system, the counting measure on the
$l$-fibres.\fidi

\smallskip
The space of units of $\G_F^{bs}(\C P_n,\pi_t)$ is the following subset of the standard simplex
$$
\Delta_n^\Z(t) = \{ c\in\Delta_n\ | \ c_k=1-t \pm e^{-\hbar n_k}\ ,\; n_k\in\Z\}  \; .
$$

\begin{figure}[!ht]
  \centering
    \includegraphics[width=1.\textwidth]{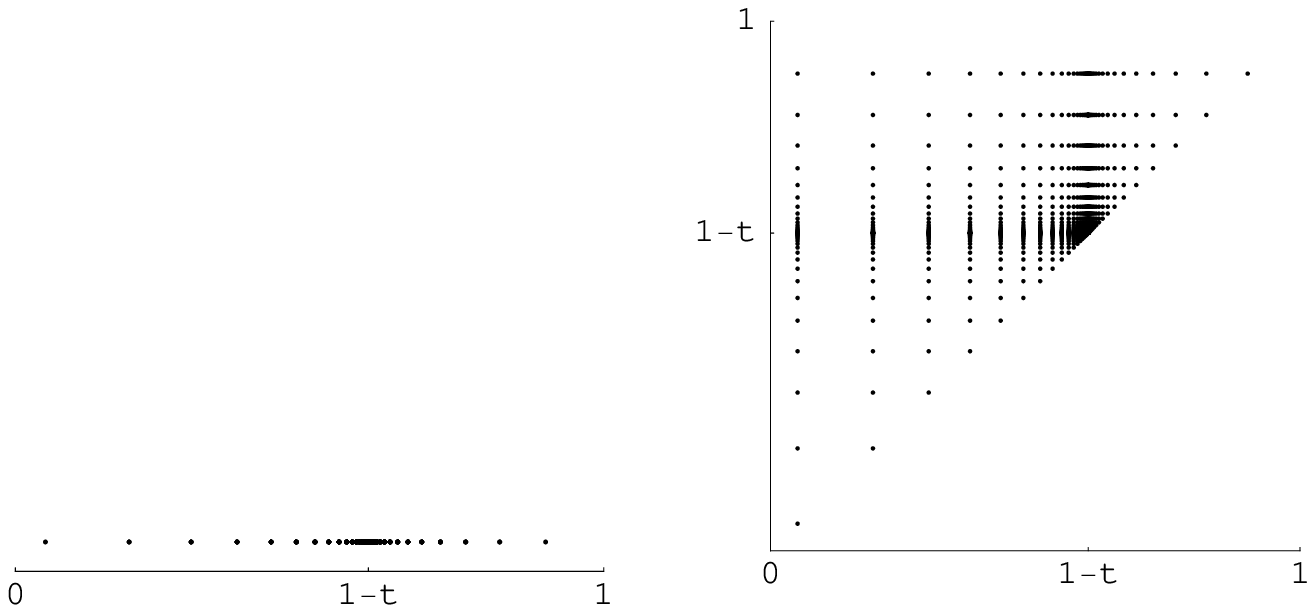}
    \caption{$\Delta_1^\Z(t)$ and $\Delta_2^\Z(t)$}
\end{figure}

For $\hbar >0$, the quasi invariant measure on $\Delta_n^\Z(t)$ associated to $f_{FS}$ is computed as
$$
\mu_{FS}(c) = \det (J_\lambda(c) + t)=\Pi_k |c_k-1+t| = \exp \left(-\hbar\sum_{k=1}^n n_k\right)  ~.
$$

\smallskip
The components of $\G_F^{bs}(\C P_n,\pi_t)$ are labeled by the couples $(r,s)$ of non negative integers such
that $r+s\leq n$. Indeed the component $\G_F^{bs}(\C P_n,\pi_t)^{(r,s)}$ is the restriction of $\G_F^{bs}(\C P_n,\pi_t)$ to
$c_k=1-t$ for $r< k\leq n-s$. Let us describe this component more explicitly.
Let $\Z^n$ act by translation on ${\overline\Z}^n$, where
$\overline{\Z}=\Z\cup\{\infty\}$ be the one point compactification. Let $\overline{\Z}^n\rtimes\Z^n$ be the action groupoid. Let
$\Delta_{r,s}^\Z(t)\subset \overline{\Z}^n$ be defined as
\begin{eqnarray}
\{(m,\infty\ldots\infty,n) \in\overline{\Z}^r\times\overline{\Z}^{n-r-s}\times\overline\Z^s|
& &  -\frac{1}{\hbar}\log (1-t) \leq m_i\leq m_{i+1}, \cr
& &   n_i\geq n_{i+1}\geq -\frac{1}{\hbar}\log t\} ~~~.
\end{eqnarray}

Note that $\Delta_{r,s}^\Z(t)$ corresponds to the stratum $\C P_n^{(r,s)}$ of dimension $2(r+s)$ of $(\C P_n,\pi_t)$ described in
Proposition \ref{foliation_pi_t} for $t\in (0,1)$. Indeed one can check that (\ref{rs_stratum}) can be restated as
$$
c_r< 1-t=c_{r+1}=c_{r+2}=\ldots = c_{n-s}< c_{n-s+1}~~~,
$$
whose $BS$ leaves are exactly those in $\Delta_{r,s}^\Z(t)$. Analogously $\Delta_{r,0}(0)$ corresponds to the $2r$-dimensional Bruhat cell for $t=0$.

For $t\in (0,1)$, the component $\G_F^{bs}(\C P_n,\pi_t)^{(r,s)}$ is isomorphic to the following subgroupoid of the action groupoid
$\overline{\Z}^n\rtimes\Z^n|_{\Delta_{r,s}^\Z}$
$$\{(q,p)\in \overline{\Z}^n\times\Z^n|_{\Delta_{r,s}^\Z}\ | ~ q_i=q_{i+1}=\infty \implies p_i=p_{i+1}\}.$$
If $t=0$ there are only components like $\G_F^{bs}(\C P_n,\pi_t)^{(r,0)}=\G^{bs,(r,0)}$ that is isomorphic to the following subgroupoid of the action groupoid
$\overline{\Z}^n\rtimes\Z^n|_{\Delta_{r,0}^\Z}$
$$\{(q,p)\in \overline{\Z}^n\times\Z^n|_{\Delta_{r,0}^\Z}\ | ~ q_i=\infty \implies p_i=0\}.$$
The description for $t=1$ is analogous.

\bigskip

\subsection{Quantization of Poisson morphisms}\label{qmorphism}

We discuss in this Section the covariance of the correspondence that sends $(\C P_n,\pi_t)$ to the \'etale groupoid $\G_F^{bs}(\C P_n,\pi_t)$,
established in Theorem \ref{main_theorem}, with respect to some of the Poisson isomorphisms that we discussed in previous Sections.

Let us first define the map $\psi:\overline{\Z}^n\rightarrow\overline{\Z}^n$ as $\psi(q)_i=q_{n+1-i}$.
The following Lemma, proved by a direct check, shows that the Poisson isomorphism between $(\C P_n,\pi_t)$ and $(\C P_n,-\pi_{1-t})$ in
Lemma \ref{tinminust} is quantized by a groupoid isomorphism.

\begin{lemma}
The map from $\in\G_F^{bs}(\C P_n,\pi_t)^{(r,s)}$ to $\G_F^{bs}(\C P_n,\pi_{1-t})^{(s,r)}$ sending $(q,p)$ to $(\psi(q),\psi(p))$ extends to a
groupoid isomorphism $\hat{\psi}: \G_F^{bs}(\C P_n,\pi_t)\rightarrow\G_F^{bs}(\C P_n,-\pi_{1-t})$ for each $t\in[0,1]$.
\end{lemma}

We are then going to show that the Poisson embeddings
$$i_k: P_k(t)\equiv p_t(S(U(k)\times U(n+1-k))\rightarrow (\C P_n,\pi_t)$$
described in Section
\ref{Poisson_geometry_complex_projective_spaces} are quantized to subgroupoids of $\G_F^{bs}(\C P_n,\pi_t)$. As a result we get at the same time the
quantization of these Poisson structures and show the covariance of the quantization with respect to these Poisson morphisms.

We know from general results that $\G(\C P_n,\pi_t)|_{P_k(t)}=\{[g\gamma] ~|~ c_k[g] = 1-t\}$ is coisotropic and the symplectic groupoid
integrating $P_k(t)$ is obtained as coisotropic reduction. The coisotropic foliation is generated by the hamiltonian flux of $c_k$
$$
\G(P_k(t),\pi_t) = \G(\C P_n,\pi_t)// \{c_k=1-t\}~~~~.
$$

The set of hamiltonians $F=\{c_i, i\not=k, h_j\}$ clearly descends to a set defining a multiplicative integrable model. All the analysis of the previous
Section descends to the quotient, so that we get the following Proposition.

\smallskip
\begin{proposition}
\label{quantization_poisson_submanifolds}
For $t\in[0,1]$, the Poisson submanifold $P_k(t)$ is quantized by the subgroupoid
$$
\G_F^{bs}(P_k(t),\pi_t) = \left\{(c,h)\in \G^{bs}_F(\C P_n,\pi_t)~|~ c_k=1-t\right\}~~~~~.
$$
\end{proposition}

\medskip
By applying the contravariant functor assigning to a groupoid its $C^*$-algebra, we get
a family of $C^*$--algebra short exact sequences
\[
0\to {\cal I}_k\to C^*(\G_F^{bs}(\C P_n,\pi_t))\to C^*(\G_F^{bs}(P_k(t),\pi_t)) \to 0~~~.
\]

It can be remarked that putting $c_{n-1}=1-t$ in $\G_F(P_n(t),\pi_t)$, and then $c_{n-2}$ and so on, we obtain a chain of subgroupoids,
that corresponds to the chain of Poisson embeddings of spheres appearing as the right side of (\ref{strataPoisson}).

\bigskip
\bigskip

\section{Quantum homogeneous spaces as groupoid \\$C^*$-algebras}

In this Section we will discuss how the quantization of the Poisson pencil and Poisson submanifolds that we described in Theorem \ref{main_theorem}
and Proposition \ref{quantization_poisson_submanifolds} is related to quantum spaces appearing in quantum group theory.

The family of quantum homogeneous spaces $\C P_{n,q,t}$ was first introduced in \cite{DiNo}, generalizing previous work by Podle\'s on $\C P_1$.
In such paper a parametrized family of right coideal $*$--subalgebras of the Hopf algebra ${\cal A}(SU_q(n+1))$ of representative functions was
introduced. The algebra ${\cal A}(\C P_{n,q,t})$ is then defined as the subalgebra of coinvariants and is characterized in terms of
generators and relations: it is a noncommutative deformation of the algebra of polynomial functions on $\C P_n$ whose semiclassical limit gives the Poisson
structure $(\C P_n,\pi_t)$.
In a similar manner the odd-dimensional quantum spheres ${\cal A}(\bbs^{2n-1}_q)$ had already been previously defined as right coideal $*$--subalgebras of  
${\cal A}(SU_q(n+1))$ in \cite{VaSo} as one of the first examples of quantum homogeneous spaces. It has to be remarked here that quantum odd spheres and standard 
complex projective spaces are quotients by quantum subgroups, while in the non standard case the quotient procedure requires the use of coisotropic quantum subgroups, 
in perfect agreement with what happens at the semiclassical level.

There is a canonical way to complete both these polynomial algebras to the $C^*$--algebra $C(\C P_{n,q,t})$ and ${\cal C}(\bbs^{2n-1}_q)$. The $C^*$--norm is defined as
$\| a\|=\mathrm{sup}_\rho\| \rho(a) \|$ of the operator norms among bounded $*$--representations $\rho$ of the polynomial algebra.

In \cite{Sheu1} Sheu, relying on the representation theory of ${\cal A}(\C P_{n,q,t})$, proved that such algebra is in fact, the groupoid
$C^*$--algebra of a topological groupoid when $t=0,1$, and also, for arbitrary $t$ in the case $n=1$. In the more general case in which $n>1$ and
$0<t<1$ his results were less conclusive though, basing on the same construction, he was able to show that $C(\C P_{n,q,t})$ is contained in a
groupoid $C^*$--algebra. On the other hand he was able to prove the existence of a short exact sequence of $C^*$--algebras 
$0\to {\cal I}\to C(\C P_{n,q,t})\to C(\bbs ^{2n-1}_q)\to 0$ relating the non standard quantum complex projective spaces to odd--dimensional quantum spheres.

We are going to show how our quantization of $(\C P_n,\pi_t)$, for $t\in[0,1]$, is related to these results. In the cases covered by Sheu's
description like $(\C P_n,\pi_0)$ and $(P_k(t)=\sphere^{2n-1},\pi_t)$,
for $k=1,n$, we will then establish that our quantization of the symplectic groupoid produces the same $C^*$-algebra of the quantum spaces obtained
via quantum group techniques. We conjecture that this is also true for $(\C P_n,\pi_t)$ and
$(P_k(t),\pi_t)$, for $t=(0,1)$ and $1<k<n$.

\subsection{The standard $\C P_{n,q,0}$}

In \cite{Sheu1}, it was shown that the $C^*$-algebra of standard quantum projective spaces is the $C^*$-algebra of a groupoid ${\mathfrak T}_n$ defined as follows.
Let us consider the translation action of $\Z^n$ on ${\overline \Z}^n$ (trivial on the points at infinity) and restrict the action groupoid to
${\overline \N}^n$:
$$
{\cal T}_n=\Z^n\times{\overline \Z}^n|_{{\overline \N}^n}=\{(j;k)\in \Z^n\times{\overline \Z}^n\,\big|\, k_i,j_i+k_i\ge 0\}\;.
$$
Consider then the subgroupoid $\tilde{\cal T}_n\subset{\cal T}_n$ given by:
$$
\tilde{\cal T}_n =\{(j;k)\in {\cal T}_n |~ k_i=\infty \implies \sum_{k=1}^i j_k=0,~ j_l=0, ~l\geq i+1\}~~~~.
$$
Finally let ${\mathfrak T}_n$ be the quotient $\tilde{\cal T}_n/\sim$ where $(j;k)\sim(j;k_1,\ldots,k_{i-1},\infty,\ldots,\infty)$ if $k_i=\infty$.
We denote with $(j;k)^{\sim}$ the equivalence class.
It is shown in Theorem 1 of \cite{Sheu1} that $C(\C P_{n,q,0})\simeq C^*({\mathfrak T}_n)$. The choice of this particular description of the groupoid
$C^*$--algebra depends on a
representation-theoretic  construction of $C(\C P_{n,q,0})$ which could be seen as an analogue of the
symplectic foliation description in terms of Lu's coordinates (see Appendix \ref{Review_Lu_coordinates}).
\begin{proposition}\label{isomorphism_standard_projective}
The map $\Phi: \G_F^{bs}(\C P_n,\pi_0)\rightarrow {\mathfrak T}_n$, defined as
$$
\Phi(c,h)=(p_1,p_2-p_1,\ldots;q_1,q_2-q_1,\ldots)^\sim, ~~~~ c_i=1-e^{-\hbar q_i},~ h_i = \hbar p_i,
$$
is an isomorphism of topological groupoids.
\end{proposition}

{\it Proof}. The map is easily seen to be well defined once one recalls the range of $(c,h)$ described in
Proposition \ref{multiplicative_integrable_model} and the characterization of Bohr-Sommerfeld leaves given in
Theorem \ref{main_theorem}. Remark that in the definition of the map $\Phi$ there is no need of specifying
the definition of $q_i-q_{i-1}=\infty-\infty$ for $i\geq r+1$, when $c_{r-1}<c_r=1$. Indeed, since $q_{r}-q_{r-1}$ is well defined and equal to
$\infty$, these entries are not relevant in the equivalence class $\Phi(c,h)$.

It is an easy check to verify that $\Phi$ is a continuous groupoid morphism. The inverse is given by
$$
\Phi^{-1}(j;k)^\sim = (c,h),~~~~ c_i=1-e^{-\hbar \sum_{l=1}^ik_l},~ h_i = \hbar \sum_{l=1}^i j_l,
$$
that is well defined and continuous.
\fidi

Let us also remark that the functorial quantization of Poisson submanifolds $P_1(0),\ldots,P_n(0)$, in this case, exactly gives rise to the
composition sequence described in Corollary 3 of \cite{Sheu1}, as a simple checking shows.

\medskip

\subsection{The non standard $\C P_{1,q,t}$, for $t\in(0,1)$}

Let us now move to the non standard case. In such case Sheu produced an explicit groupoid description only for $n=1$. In the $n>1$ case, the non
standard quantum complex projective space is described just as a sub $C^*$--algebra of a groupoid $C^*$--algebra.

For $n=1$ in \cite{Sheu1} it was proved that $C(\C P_{1,q,t})\simeq C(\mathfrak G)$, $t\in(0,1)$, where
$\mathfrak G$ is the subgroupoid of the action groupoid $\Z^2\times \overline{\Z}^2|_{\overline \N^2}$
\[
{\mathfrak G}=\{(j,j,k_1,k_2)\in \Z^2\times \overline{\Z}^2|_{\overline \N^2}\,|\, k_1~ {\rm or} ~k_2=\infty\}\, .
\]
Let us recall that ${\cal G}^{bs}_{F}(\C P_1,\pi_t)$ has three components
\[
\begin{array}{rcl}
{\cal G}^{bs}_{F}(\C P_1,\pi_t)^{(1,0)}&=&\{(q,p)\in {\overline \Z}\times \Z\,|\, q\ge q_t,\, q+p\ge q_t\,\}\\
{\cal G}^{bs}_{F}(\C P_1,\pi_t)^{(0,1)}&=&\{(q,p)\in {\overline \Z}\times \Z\,|\, q\ge {\tilde q}_t,\, q+p\ge {\tilde q}_t\,\}\\
{\cal G}^{bs}_{F}(\C P_1,\pi_t)^{(0,0)}&=&\{(\infty,p)\in {\overline \Z}\times \Z\,\}
\end{array}
\]
where $q_t=-\frac 1 \hbar \ln(1-t)$ and ${\tilde q}_t=-\frac 1 \hbar \ln t$. A direct computation
shows that:
\begin{proposition}
The map $\Phi:\G_F^{bs}(\C P_1,\pi_t)\rightarrow {\mathfrak G}$ defined on each component by:
\[
\Phi(q,p)=\left\{
\begin{array}{ccc}
&(p,p,\infty, q-q_t)& \, \mathrm{ if}\, (q,p)\in {\cal G}^{bs}_{F}(\C P_1,\pi_t)^{(1,0)}\\
&(p,p,q-{\tilde q}_t,\infty)& \, \mathrm{ if}\, (q,p)\in {\cal G}^{bs}_{F}(\C P_1,\pi_t)^{(0,1)}\\
&(p,p,\infty, \infty)& \, \mathrm{ if}\, (q,p)\in {\cal G}^{bs}_{F}(\C P_1,\pi_t)^{(0,0)}
\end{array}
\right.
\]
is an isomorphism of topological groupoids for all $t\in (0,1)$.
\end{proposition}

The existence of this isomorphism also proves that the groupoids of BS-leaves and so the groupoid $C^*$--algebras are not depending on $t$, as long as
$0<t<1$.

\medskip

\subsection{Odd dimensional quantum spheres $\sphere^{2n-1}_q$}\label{quantum_spheres}

From the discussion in Section \ref{Poisson_geometry_complex_projective_spaces} we know that the Poisson submanifold $P_n(t)\subset (\C P_n,\pi_t)$
is isomorphic
to the odd Poisson sphere $\sphere^{2n-1}$. We want
to compare the groupoid $\G_F^{bs}(P_n(t),\pi_t)$ appearing in Proposition \ref{quantization_poisson_submanifolds} with the characterization
of the quantum sphere given
in \cite{Sh-sphere}, where it is shown that $C(\sphere^{2n-1}_q)\simeq C^*({\mathfrak F}_{n-1})$, where
${\mathfrak F}_{n-1}$ is a groupoid that we are going to describe.

Let ${\cal F}_{n-1}$ be the following subgroupoid of the action groupoid (with the first $\Z$ acting trivially)
\[
\{
(z;x;w)\in \Z\times (\Z^{n-1}\times {\overline \Z}^{n-1})|_{\overline \N^{n-1}}| \,
w_i=\infty\Rightarrow z=-\sum_{j} x_j, \ x_{k}=0, k>i
\}
\]
Define on it the equivalence relation  $(z;x;w)\sim (z;x;w_1,\ldots,w_i,\infty)$ if $w_i=\infty$  and let $\mathfrak F_{n-1}={\cal F}_{n-1}/\sim$.
Let us denote with $(z;x;w)^\sim\in{\mathfrak F}_{n-1}$ the equivalence class. Let, as before, $q_t=-\frac{1}{\hbar}\log(1-t)$.

\begin{proposition}
\label{isomorphism_odd_spheres}
The following map $\Phi:\G^{bs}_F(P_n(t),\pi_t)\rightarrow {\mathfrak{F}_{n-1}}$ defined as follows
$$
\Phi(c,h) = (-p_n; p_1,p_2-p_1,\ldots,p_{n-1}-p_{n-2}; q_1-q_t, q_2-q_1,\ldots,q_{n-1}-q_{n-2})^\sim~,
$$
where $h_k = \hbar p_k$ and $c_k=1-t-e^{-\hbar q_k}$ for $k<n$, is an isomorphism of topological groupoids.
\end{proposition}
{\it Proof}. Let us recall that if $(c,h)\in \G_F^{bs}(P_n(t),\pi_t)$ then $c_n=1-t$. If $c_{i-1}<c_i=1-t$, for some $i$ then
$w_i=q_i-q_{i-1}=\infty$ is well defined and is not relevant how we define
$w_r=q_r-q_{r-1}=\infty-\infty$ for $r>i$ because of the equivalence relation. Furthermore, in this case $x_k=p_k-p_{k-1}=0$ for $k>i$ and
$\sum_j x_k = p_i=p_{i+1}=\ldots = p_n=-z$.
The inverse $\Phi^{-1}:{\mathfrak{F}_{n-1}}\rightarrow \G_F^{bs}(\C P_n,\pi_t)$ is given by
$$
c_i=(1-t)(1- e^{-\hbar\sum_{j=1}^iw_j }),~ h_i=\hbar \sum_{j=1}^i x_j,~ c_n=1-t,~  h_n=-\hbar z~~~.
$$
The proof that $\Phi$ is a groupoid morphism is straightforward. \fidi

\medskip
Since the functor assigning to a groupoid its $C^*$--algebra is contravariant there is a $C^*$--algebra epimorphism from the
$C^*$--algebra of $\C P_{n,q,t}$, for $t\in(0,1)$, to the $C^*$--algebra of the odd-dimensional quantum sphere $\sphere^{2n-1}_q$. This is exactly
the analogue of Theorem 6 in
\cite{Sheu1}.

\bigskip
\bigskip
\section{Conclusions}

\noindent{\it Quantum homogeneous spaces as groupoid $C^*$-algebras}. In Theorem \ref{main_theorem} we proved that the groupoid of Bohr-Sommerfeld
leaves $\G_F^{bs}(\C P_n,\pi_t)$ has a unique Haar system for all $t\in[0,1]$. According to our program, we can define the convolution algebra and
the $C^*$-algebra out of these data. For $t=0$, by using Sheu's results we know that this is exactly the $C^*$-algebra defined from the quantum
homogeneous space $\C P_{n,q,0}$. In the case $t\in (0,1)$ Sheu identifies the quantum algebra only as a subalgebra of a groupoid $C^*$-algebra.
It would be very interesting to have a proof based on the representation theory of $\C P_{n,q,t}$ that its $C^*$-algebra $C(\C P_{n,q,t})$ is isomorphic to
$C^*(\G_{F}^{bs}(\C P_n,\pi_t))$ also for $t\in(0,1)$. Analogously it would be interesting to determine the polynomial algebra of quantum $P_k(t)$ and show
that its $C^*$-algebra is isomorphic to $C^*(\G_F(P_k(t),\pi_t))$ (we discussed the case $k=1,n$ in Subsection \ref{quantum_spheres}).

\smallskip

\noindent{\it Compact Hermitian symmetric spaces}. The basic ingredient of the construction that we have presented in this paper is the bihamiltonian system on $\C P_n$. 
This system exists for all compact hermitian symmetric spaces; most of the results extend to the general case. In particular, the proof of the multiplicativity of the 
modular function is general; given the eigenvalues of the Nijenhuis operator, the structure maps of the groupoid of leaves are expressed by the same formulas.

\smallskip

While the existence of the bihamiltonian system is known since long time, at the best of our knowledge, $\C P_n$ is the only case that has
been investigated. The key characterization that we need is the diagonalization of the Nijenhuis operator: for $\C P_n$ is a result easy to obtain
(Proposition \ref{action_poisson_pencil}). In the general case we don't know if global eigenvalues of the Nijenhuis operator exist, how to describe them and
finally if the whole quantization procedure works. This requires further investigation, that makes a reasonable plan for our future research.

\smallskip

\noindent{\it Prequantization cocycle}. The prequantization cocycle did not appear in our final quantization setting; this means that we produce
our quantum algebras by taking it to be trivial. On one side, in
those cases where we can compare our results with quantum group constructions, via Sheu's results, there is no need of a non trivial 2-cocycle. On the other
side we know that, at least in some examples but we expect it in general, its cohomology class is not trivial. It is possible that it makes a
difference what kind of cohomology we consider, differentiable or just continuous, and that for the $C^*$-algebras smooth properties are not relevant.
An indication in this direction is given by recent results of \cite{Sheu4} where it is shown that certain Poisson-homogeneous structures which are
cohomologically non trivial may become trivial if one admits vector fields which are smooth along the symplectic leaves and globally only continuous.
In any case, the role of this prequantization cocycle is still unclear.

\smallskip

\noindent{\it Morita equivalences}. The Poisson structures $(\C P_n,\pi_t)$ are not in general Poisson diffeomorphic for different $t$.
In \cite{BR,Rad} it was proven that non standard $\C P_1$ Poisson structures are non diffeomorphic
for different values of $t$. In fact the Poisson invariant allowing to distinguish them is the so called regolarized Liouville volume
$V(\pi)=\ln\frac{1-t}{1+t}$. However, since the modular periods are the same, apart from the singular $t=0,1$ case,  non standard Poisson
$\C P_1$ are all Poisson-Morita equivalent. The fact that the corresponding groupoid of Bohr-Sommerfeld leaves does not depend, as topological groupoid, on $t$ can be seen 
as a quantum  analogue of this fact.
It would be interesting to analyze these properties for higher $n$.  We remark that the example developed in \cite{BCT} shows that even non
Poisson-Morita equivalent structures can produce the same topological groupoid of BS leaves.
We expect that, like in the $n=1$ case, all $(\C P_n,\pi_t)$ are at least Morita equivalent. We plan to investigate this equivalence and how it is
quantized to the Bohr-Sommerfeld groupoids.

\bigskip

\noindent{\bf Acknowledgments}
The research of J.Q. is supported by the Luxembourg FNR grant PDR 2011-2, and by the UL Grant GeoAlgPhys 2011-2013
\bigskip

\appendix

\section{Review of Lu's Coordinates}\label{Review_Lu_coordinates}
In this section, we recall the derivation of Lu's parametrization of Bruhat cells \cite{Lu1}. \emph{Do notice that, in contrast to the main text, the homogeneous spaces 
such as the flag manifolds are presented as right cosets, e.g. $\BB{C}P_n=SU(n+1)/U(n)$.} This is to avoid certain awkwardness in the presentation, no conclusions in the main 
text will be changed.

We focus on the case
$K=SU(n)$, $K^*=SB(n)$ and $G=K_{\BB{C}}=SL(n,\BB{C})$, for such groups, the concepts such as the Bruhat Poisson structure, Bruhat cells, dressing actions,
etc can all be understood rather concretely.

Our focus is the flag manifold flag$(n_1,\cdots,n_k)$, where the string of increasing integers $0=n_0< n_1<n_2<\cdots<n_{k}<n_{k+1}=n$ denotes the nested subspaces of $\BB{C}^n$ of dimension
$n_1,n_2,\cdots,n_k$. A point in the flag can be represented by a matrix $g\in SL(n)$, suppose $g=[v_1,\cdots v_n]$, where $v_i$ are column vectors, then the
flag in question is $\textrm{span}(v_1,\cdots v_{n_1})\subset\textrm{span}(v_1,\cdots v_{n_2})\subset \textrm{span}(v_1,\cdots v_{n_{k}})\subset\BB{C}^n$.
There are certain redundancies of such a parametrization, for example, a right multiplication by a matrix in $SB$
does not change the flag, since these matrices are upper-triangular. In general, unless $n_i-n_{i-1}=1$ for all $1\leq i\leq k$, i.e. the full-flag manifold, the stability group is bigger
 than $SB$.

In contrast, the left multiplication by $\gc\in SB$ on $g=[v_1,\cdots v_n]\in SU$ induces a non-trivial action on the flag, so one can decompose the flag manifold into orbits of this action.
 The action is related to the dressing action reviewed in the text, write $\gc\cdot g={}^{\gc}g\cdot \gc^g$, with ${}^{\gc}g\in SU$ and $\gc^g\in SB$. One can discard $\gc^g$ since
${}^{\gc}g\cdot \gc^g$ and ${}^{\gc}g$ give the same flag. If the stability group were exactly $SB$, then the action can be identified with the dressing transformation.

To cut down the stability group one can do the following. The identity matrix corresponds to $[e_1,\cdots,e_n]$, where $e_i$ are the standard basis vectors, while an element of the Weyl
group $W$, when represented as a matrix, corresponds to $[e_{\gs(1)},\cdots,e_{\gs(n)}]$, where $\gs$ is a permutation (recall that the Weyl group for $SU$ are none other than permutations).
 One can certainly reach any point in the flag by left multiplying $SB$ to a suitable $w\in W$, thus one can label the orbits of the left $SB$ action by elements of $W$. Further introduce a
 subset of permutations called shuffles, denoted as $\gs_n(n_1,n_2,\cdots,n_k)$, consisting of $\{\gs\in S_n\,|\,\gs(i)<\gs(j)~\textrm{if}~n_k<i<j\leq n_{k+1}\}$. It is clear that the
 shuffles $\gs(n_1,\cdots,n_k)$ label all the orbits of left $SB$ action in the flag, since permutations within each subspace spanned by, say, $\{e_{n_i+1},\cdots e_{n_{i+1}}\}$ does not
change the flag. The point is that now one may impose a 'gauge condition' on the
representative of a point in the flag, as follows. One defines an ordering on the vectors $\{v_i\}$ that span the flag by saying that $v<v'$ if the last nonzero component of $v'$ appears
 later than that of $v$. The gauge condition is such that: if $n_k< i<j\leq n_{k+1}$, then $v_i<v_j$. It is clear that the image of the left $SB$ action on a shuffle satisfies this gauge
 condition. Once this condition is enforced, the stability group is cut down to being the right multiplication by $SB$. To summarize, the flag has the following decomposition
\bea \textrm{flag}(n_1\cdots,n_k)=\bigcup_{w}~SB\,w\,SB/SB,\nn\eea where the summation is over $\gs_n(n_1,n_2,\cdots,n_{k})$-shuffles.

Next, we shall identify a subgroup $N_w\subset SB$ whose left action on a particular shuffle $w$ is free. Let $\xi\in\FR{sb}$, then the action of $\xi$ on $w$ is nonzero if
$p_1Ad_{w^{-1}}\xi\neq0$, where $p_1:\;\FR{sl}\to\FR{su}$ is the projection. Recall that $SB$ is a subgroup generated by the positive roots of $\FR{su}$ and $i\FR{t}$, where $\FR{t}$ is
the Cartan. One sees that for a given $w$, the subset of roots $\Phi_w=\{\ga>0|Ad_{w^{-1}}\ga<0\}$ will generate a nonzero action since $p_1Ad_{w^{-1}}\ga\neq0$, and thus $N_w$ is
generated by $\Phi_w$.

One can always decompose an element $w$ of the Weyl group into a string of Weyl reflections $w=w_{\ga_k}\cdots w_{\ga_1}$, where each $\ga_i$ is a simple root and $w_{\ga_i}$ is the Weyl
 reflection of $\ga_i$. If we denote by $C_w$ the image of the left $SB$ action on $w$, then Lu showed that there is a diffeomorphism
\bea C_{\ga_k}\times\cdots \times C_{\ga_1}\to C_w\label{orbit_prod},\eea
where we have abbreviated $C_{w_{\ga_i}}$ as $C_{\ga_i}$. The map is given by taking the product of the lhs, furthermore, this map is a Poisson diffeomorphism, where the lhs is given
the product Poisson structure. This gives us a very practical way of parametrizing the orbit $C_w$ by that of $C_{\ga}$. Recall that a simple root $\ga$ gives an embedding of $SU(2)$
into $SU(n)$, so we may focus on the $SU(2)$ case.

For $SU(2)$, the non-trivial Weyl group element can be chosen as $w=i\gs_2$, where $\gs_i$ are  Pauli matrices. There is only one simple root $\gs_+=(\gs_1+i\gs_2)/2$, and
$w^{-1}\gs_+w=-\gs_-$, so $N_w$ is generated by $\gs_+$, and hence can be parametrized as
\bea \gc(y)=e^{y\gs_+}=\left(
       \begin{array}{cc}
         1 & y \\
         0 & 1 \\
       \end{array}
     \right),~~~y\in\BB{C}\nn\eea
it is then easy to work out the orbit of the dressing action
\bea \gc \,w=\epsilon^{-1}\left[\begin{array}{cc}
       -y & 1 \\
       -1 & -\bar y
     \end{array}\right]\left[\begin{array}{cc}
       \epsilon & 0 \\
       0 & \epsilon^{-1}
     \end{array}\right]\left[\begin{array}{cc}
       1 & -\epsilon^{-2}\bar y \\
       0 & 1
     \end{array}\right],~~~~\epsilon=(1+|y|^2)^{1/2},\nn\eea
where the last two factors belong to $SB$. Thus the orbit of the dressing action is parametrized as
\bea g(y)=\frac{1}{\epsilon}\left[\begin{array}{cc}
       -y & 1 \\
       -1 & -\bar y
     \end{array}\right],~~~~\epsilon=(1+|y|^2)^{1/2}.\label{SU(2)_cell}\eea
This orbit actually covers all points of the flag $\BB{C}P_1$ except the point $z=\infty$, i.e. the unique dimension 2 Schubert cell of $\BB{C}P_1$.

To work out the general formula Eq.\ref{orbit_prod} for $\BB{C}P_n$, we fix the simple roots for $SU(n+1)$ to be: $\ga_i=[0,\cdots \underbrace{1,-1}_{i,i+1},0, \cdots 0]$, then the Weyl
 reflection $w_{\ga_i}$ is an identity matrix except that the $(i,i+1)$ diagonal block is occupied by $i\gs_2$.
Next we work out as an example the parametrization of the largest cell of $\BB{C}P_n$. The Weyl group element corresponding to this cell is the shuffle $\gs_{n+1}(1)$, and it can be
 decomposed into
\bea \gs_{n+1}(1)=w_{\ga_n}\cdots w_{\ga_{1}}.\label{decomp_shuffle}\eea
The orbit of each $C_{\ga_i}$ is parametrized as a rank $(n+1)$ block diagonal matrix,
\bea {\tt g}_i(y_i)=\left[\begin{array}{ccc}
       {\bf 1}_{i-1} & 0 & 0 \\
       0 & g(y_i) & 0\\
       0 & 0 & {\bf 1}_{n-i} \\     \end{array}\right],\label{g_i}\eea
where $g(y)$ is defined in Eq.\ref{SU(2)_cell}.
It is then not hard to compute the product Eq.\ref{orbit_prod} explicitly, in fact, we will only need the first column of $C_w$, which is given by
\bea&& \left[{\tt g}_n(y_n)\cdots {\tt g}_1(y_1)\right]_{\sbullet,1}\nn\\
&
=&\frac{1}{\epsilon_1\cdots \epsilon_n}\left[
       \begin{array}{cccccc}
         (-1)^{n}y_1\epsilon_2\cdots\epsilon_n, & \cdots & -y_{n-2}\epsilon_{n-1}\epsilon_n, & y_{n-1}\epsilon_n, & -y_n, & -1 \\
       \end{array}\right]^T,\nn\eea
where $\epsilon_i=(1+|y_i|)^{1/2}$. This first column of $C_w$ should be identified with the homogeneous coordinates $[X_1,\cdots X_{n+1}]$ of $\BB{C}P_n$, and $C_w$ corresponds to
the cell where $X_{n+1}\neq0$. For convenience of the identification, we will make a redefinition of variables $y_{n-1}\to -y_{n-1},\;y_{n-3}\to -y_{n-3},\cdots$, then we have
\bea z_n=\frac{X_n}{X_{n+1}}=y_n,~~z_{n-1}=\frac{X_{n-1}}{X_{n+1}}=y_{n-1}\epsilon_n,\cdots.\nn\eea
And we have the inverse map
\bea y_n&=&z_n,\nn\\
y_{n-1}&=&z_{n-1}(1+|z_n|^2)^{-1/2},\nn\\
y_{n-2}&=&z_{n-2}(1+|z_{n-1}|^2+|z_n|^2)^{-1/2},\nn\\
&\cdots&\nn\\
y_1&=&z_1(1+|z_2|^2+\cdots+|z_n|^2)^{-1/2}.\label{Lu_coord}\eea
These are the Lu's coordinates quoted in \eqref{lu_coordinates}. They can be thought as coordinates on leaves of different copies of $SU(2)$ projected onto the complex projective space. 
Let us remark that in principle all that is contained in here applies also to the nonstandard case, the only difference being in the projection map $p_t$, which will pick a $t$-depending 
linear combination of the first and last column of $C_w$.

\subsection{The Kirillov-Kostant symplectic form}
In the text, we shall also need the symplectic Kirillov-Kostant symplectic form for $\BB{C}P_n$, we review here the construction of this form and some concrete formulae and show that
it takes a very simple form in Lu's coordinates.

For a Lie group $K$, the coadjoint orbit associated with $\gl\in\FR{g}^*$ is ${\cal O}_{\gl}=\{(Ad_g^{-1})^*\gl\,|\,g\in K\}$, and it is
diffeomorphic to $K/H_{\gl}$ where $H_{\gl}$
is the stability group of $\gl$. Over ${\cal O}_{\gl}$, one can write down a form
\bea \omega_{\gl}=id\Tr[g^{-1}dg\gl]=-i\Tr[(g^{-1}dg)^2\gl],\label{Kirillov_Kostant}\eea
It is non-degenerate and closed (but \emph{not} exact).

In the case of $K=SU(n+1)$, the coadjoint orbit passing through
\bea \gl=\frac{1}{n+1}\textrm{diag}[n,\underbrace{-1,\cdots,-1}_n]\label{highest_weight_cpn}\eea
is nothing but $\BB{C}P_n$, since $H_{\gl}$ in this case is clearly $S(U(1)\times U(n))$.

For $SU(2)$ the computation is simple, letting $\gl=[1/2,-1/2]$ and parametrize a group element as
$g(y)e^{i\gt\gs_3}$ (where $g(y)$ is defined in Eq.\ref{SU(2)_cell}), one gets
\bea \omega_{\gl}=\frac{i}{2}d\Tr[g^{-1}dg\sigma_3]=\frac{-idy\wedge d\bar y}{(1+|y|^2)^2}\label{KK_SU(2)},\eea
the Fubini-Study K\"ahler form on $\BB{C}P_1$.

To obtain the Kirillov-Kostant for $\BB{C}P_n$, we take $\gl$ as in Eq.\ref{highest_weight_cpn}, i.e. the highest weight of the fundamental representation of $SU(n+1)$.
Then we shall use the Lu parametrization for $SU(n+1)/S(U(1)\times U(n))$.
The largest cell in this quotient is labeled by the shuffle $w=\gs_{n+1}(1)$, and it can be decomposed as in Eq.\ref{decomp_shuffle}
The cell $C_w$ is now a product
\bea C_w=C_{\ga_n}\cdots C_{\ga_2}C_{\ga_1},\nn\eea
We parametrize this product $C_{\ga_n}\cdots C_{\ga_2}C_{\ga_1}$
as ${\tt g}_n(y_n)\cdots {\tt g}_1(y_1)$, where ${\tt g}_i(y_i)$ is defined in Eq.\ref{g_i}.

One may now proceed to compute the Liouville form $\Tr[g^{-1}dg\gl]$ leading to
\bea &&\Tr[g^{-1}dg\gl]=iA(z_1)+\frac{i}{\epsilon_1^2}A(y_2)+\frac{i}{\epsilon_1^2\epsilon_2^2}A(y_3)+\cdots \frac{i}{\epsilon_1^2\cdots \epsilon_{n-1}^2}A(y_n),\label{Liouville_n}\\
&&\hspace{4cm}A(y)=\frac{i}{2}\frac{yd\bar y-yd\bar y}{(1+|y|^2)}\nn.\eea
Let us plug in the change of variables Eq.\ref{Lu_coord} and get the neat result
\bea \Tr[g^{-1}dg\gl]=-\frac{1}{2}\frac{\sum_{i=1}^n z_id\bar z_i -c.c}{1+\sum_{i=1}^n|z_i|^2}\nn.\eea
The differential of the above 1-form is recognized as the Fubini-Study K\"ahler form for $\BB{C}P_n$.

\subsection{Concrete formulae for the Bruhat Poisson structure}
In this section we collect some formulae for the Poisson structure of $SU(2)$ and its dual group $SB(2)$, as well as the symplectic form on the double $SL(2)$. The computation leading
to these results are quite lengthy and is omitted.

Parametrize an element of $SU(2)$ as $g(y)e^{i\gt\gs_3}$, and plug it into Eq.\ref{poisson_K}, a straightforward calculation shows
\bea \pi_{SU(2)}=i(1+|y|^2)\frac{\partial}{\partial y}\wedge\frac{\partial}{\partial\bar y}\nn,\eea
which vanishes when $y\to \infty$. As mentioned above the map Eq.\ref{orbit_prod} is a Poisson map, so the Poisson tensor for a general cell is just the sum of the results above.

The group $SB(2)$ is dual group of $SU(2)$ , and has a simple parametrization $\gc=e^{\phi\gs_3}e^{u\gs_+}$, and the dual Poisson structure in this 
coordinate system is 
\bea -2\pi_{SB(2)}=2i\big(-(1+|u|^2)+e^{-4\phi}\big)\partial_u\wedge\partial_{\bar u}+iu\partial_u\wedge\partial_{\phi}-i\bar u\partial_{\bar u}\wedge\partial_{\phi}.\nn\eea

On the double $SL(2)$, one has a Poisson structure $\Pi$ defined in Eq.\ref{Poisson_double}, which is non-degenerate, so its inverse defines a 
symplectic structure.
If one parametrizes the double as $\gc \cdot g=e^{\phi\gs_3}e^{u\gs_+}g(y)e^{i\gt\gs_3}$, then
\bea \Pi^{-1}=\frac{i}{(1+|w|^2)}dwd\bar w-\frac{i}{(1+|y|^2)}dyd\bar y-2d\gt d\log\frac{\ep^2(w)}{e^{2\phi}\ep^2(y)},\label{Pi_inverse_direct}\eea
where $w=e^{2\phi}(u+y)$. This nice structure of the symplectic form persists for a more general group, and the result and proof thereof shall appear in a forthcoming publication.

\bigskip
\bigskip

\section{An alternative proof of Lemma \ref{poisson_pencil}}\label{direct_proof}

We also present a proof using direct computation so as to be self-containing.
Denote $U_0=S(U(1)\times U(n))\subset SU(n+1)$, and $U_t=Ad_{\gs_t}U_0$. Now $U_0$ is a Poisson-Lie subgroup while $U_t$ is coisotropic. The quotient $U_t\backslash SU(n+1)$ is also 
isomorphic to $\BB{C}P_n$, and there is a well-defined map between the cosets $\phi_t:~U_t\backslash SU(n+1)\to U_0\backslash SU(n+1)$,
\bea [g]_t \stackrel{\phi_t}{\to} [\gs_t^{-1}g]_0,~~~g\in SU(n+1),\label{map_coset}\eea
where $[\;]_t$ is the quotient w.r.t $U_t$.

Since $U_t$ is coisotropic, there is a Poisson structure $\pi_t$ on $U_t\backslash SU(n+1)$ inherited from that of $SU(n+1)$. Concretely, one chooses locally a representative $g$ of 
$[g]_t$, one can then trivialize $T^*(U_t\backslash SU(n+1))$ by the right invariant 1-forms: $r_{g^{-1}}^*\xi,~\xi\in\FR{u}_t^{\perp}$, and under a change of representative 
$g\to hg,~h\in U_t$, one has $\xi\to Ad_{h^{-1}}^*\xi$. One then defines \bea \pi_t([g]_t)(r_{g^{-1}}^*\xi,r_{g^{-1}}^*\eta)=\pi(g)(r_{g^{-1}}^*\xi,r_{g^{-1}}^*\eta),\nn\eea
where $\pi$ is the Poisson tensor of $SU(n+1)$ and one can check that this definition is independent of the choice of the local representative, thanks to the multiplicativity of $\pi$.

Next, we will push $\pi_t$ forward to $U_0\backslash SU(n+1)$ using the map (\ref{map_coset}).
Let now $\xi,\eta \in\FR{u}_0^{\perp}$, and plug the 1-forms $r_{g^{-1}}^*\xi,\,r_{g^{-1}}^*\eta$ into $\phi_{t*}\pi_t$,
\bea &&(\phi_{t*}\pi_t)(g)((r_{g^{-1}})^*\xi,(r_{g^{-1}})^*\eta)=\pi(\gs_tg)((l_{\gs^{-1}_t})^*(r_{g^{-1}})^*\xi,(l_{\gs^{-1}_t})^*(r_{g^{-1}})^*\eta)\nn\\
&=&\big(l_{\gs_t*}\pi(g)+r_{g*}\pi(\gs_t)\big)((l_{\gs^{-1}_t})^*(r_{g^{-1}})^*\xi,(l_{\gs^{-1}_t})^*(r_{g^{-1}})^*\eta)\nn\\
&=&\pi(g)((r_{g^{-1}})^*\xi,(r_{g^{-1}})^*\eta)+\pi(\gs_t)((l_{\gs^{-1}_t})^*\xi,(l_{\gs_t^{-1}})^*\eta)\nn\\
&=&\pi(g)((r_{g^{-1}})^*\xi,(r_{g^{-1}})^*\eta)-\pi(\gs^{-1}_t)((r_{\gs_t})^*\xi,(r_{\gs_t})^*\eta)\label{used_later}.\eea
where in the last step we used the multiplicativity of $\pi$ and that $\pi(1)=0$.

The first term of (\ref{used_later}) gives $\pi_0$, we want to show that the second term gives $\pi_{\gl}$, namely it is the inverse of the Kirillov-Kostant symplectic form
\bea &&\omega_{\gl}=-i\Tr[(dg\cdotp g^{-1})^2\gl],~~~\gl=(n+1)^{-1}\textrm{diag}[n,\underbrace{-1,\cdots,-1}_n],\label{KK_form}\eea
note the change of convention regarding left and right invariance compared to the formula (\ref{Kirillov_Kostant}) in Appendix.\ref{Review_Lu_coordinates}.
Apply the formula (\ref{poisson_K}) for $\pi$,
\bea &&-\pi(\gs^{-1}_t)((r_{\gs_t})^*\xi,(r_{\gs_t})^*\eta)=-\frac12\bra p_1 Ad_{\gs_t}\xi,p_2Ad_{\gs_t}\eta\ket\nn.\eea
Now let
\bea \xi=\left(
               \begin{array}{cc}
                 0 & \vec v \\
                 0 & 0 \\
               \end{array}
             \right),~~~~\eta=\left(
               \begin{array}{cc}
                 0 & \vec w \\
                 0 & 0 \\
               \end{array}
             \right),~~~~~\xi,\eta\in\FR{u}_0^{\perp}\nn\eea
A direct calculation shows
\bea \bra p_1\gs_t\xi\gs_t^{-1},p_2\gs_t\eta\gs_t^{-1}\ket=-t\im(\vec v\cdot \vec w^{\dagger}),\label{earlier}\eea
thus $-\pi(\gs^{-1}_t)((r_{\gs_t})^*\xi,(r_{\gs_t})^*\eta)=t/2\im(\vec v\cdot \vec w^{\dagger})$.

Next, let $X,\,Y\in \FR{su}(n+1)$ and evaluate (\ref{KK_form}) on $r_{g*}X,\,r_{g*}Y$,
\bea \go_{\gl}(r_{g*}X,\,r_{g*}Y)=-i\Tr[[X,Y]\gl]=2\sum_{k=2}^{n+1}\im(X_{1k}Y_{k1}),\nn\eea
comparing this with (\ref{earlier}), we conclude
\bea \phi_{t*}\pi_t=\pi_0+t\omega_{\gl}^{-1}=\pi_0+t\pi_{\gl}\nn\eea
\fidi

\end{document}